# COUNTING PLANAR RANDOM WALK HOLES

By Christian Beneš

*Tufts University*

We study two variants of the notion of *holes* formed by planar simple random walk of time duration $2n$ and the areas associated with them. We prove in both cases that the number of holes of area greater than $A(n)$, where $\{A(n)\}$ is an increasing sequence, is, up to a logarithmic correction term, asymptotic to $n \cdot A(n)^{-1}$ for a range of large holes, thus confirming an observation by Mandelbrot. A consequence is that the largest hole has an area which is logarithmically asymptotic to $n$. We also discuss the different exponent of $5/3$ observed by Mandelbrot for small holes.

**1. Introduction.** The object of our study in this paper is planar simple random walk $S$ defined by $S(0) = (0,0)$ and for $n \in \mathbb{N} = \{1, 2, \ldots\}$ by $S(n) = \sum_{j=1}^{n} X_j$, where $\{X_j\}_{j \in \mathbb{N}}$ are independent random vectors satisfying $\mathbb{P}(X_j = \pm e_i) = 1/4, i = 1, 2$, where $e_1 = (1, 0)$ and $e_2 = (0, 1)$. We will also think of planar simple random walk $S$ as being a continuous process, that is, for noninteger times $t$, we let $S(t)$ be the linear interpolation of the walk's position at the surrounding integer times: For all real $t \geq 0$,

$$(1.1) \qquad S(t) = S([t]) + (t - [t])(S([t] + 1) - S([t])),$$

where $[t]$ denotes the integer part of $t$. For any real numbers $0 \leq a \leq b$, we will write $S[a,b] := \{S(t)\}_{a \leq t \leq b}$, and use the same notation for Brownian motion $B$. Let the *holes* or *components* made by $S[0, 2n]$ be the connected components of $\mathbb{C} \setminus S[0, 2n]$, where $\mathbb{C}$ denotes the complex plane, and the *lattice holes* made by $S[0, 2n]$ be the connected components of $\mathbb{Z}^2 \setminus \{S(j)\}_{j \in \{0, \ldots, 2n\}}$. Two points $z$ and $w$ lie in a same component of $\mathbb{Z}^2 \setminus \{S(j)\}_{j \in \{0, \ldots, 2n\}}$ if they can be joined by a nearest-neighbor path in $\mathbb{Z}^2$ that does not intersect $S[0, 2n]$. The *area* of a hole is defined to be its Lebesgue measure and the *lattice area* of a lattice hole is its cardinality. See Figure 1.

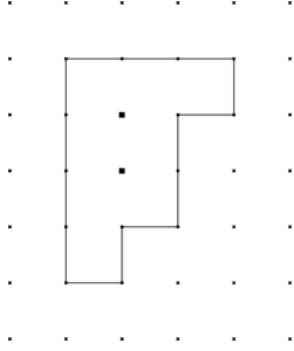

Fig. 1. *Two types of holes and area: here, area = 8, lattice area = 2.*

We define for any $r \in \mathbb{R}_+, n \in \mathbb{N}$,

$$\mathcal{H}_n(r) = \#\{\text{holes of area} \geq r \text{ made by } S[0, 2n]\},$$

(1.2) $\quad \mathcal{L}_n(r) = \#\{\text{lattice holes of lattice area} \geq r \text{ made by } S[0, 2n]\},$

$$\tilde{N}_n(\delta) = \mathcal{H}_n(n^{1-\delta}), \qquad \hat{N}_n(\delta) = \mathcal{L}_n(n^{1-\delta}).$$

Our main result is the following:

THEOREM 1.1. *There exists a $\delta_0 > 0$ such that for all $0 < \delta \leq \delta_0$,*

(1.3) $$\frac{\log^2(n^\delta)}{n^\delta} \tilde{N}_n(\delta) \xrightarrow{P} 2\pi \qquad as\ n \to \infty,$$

(1.4) $$\frac{\log^2(n^\delta)}{n^\delta} \hat{N}_n(\delta) \xrightarrow{P} 2\pi \qquad as\ n \to \infty,$$

*where $\xrightarrow{P}$ denotes convergence in probability.*

We first outline the key ideas of the argument for (1.3), which is the part of Theorem 1.1 for which we give a full proof in this paper. The proof of (1.4) is practically identical and we will just mention in Section 5 which small modifications are needed to obtain it.

1. Use Brownian scaling to extend the result in [10] (see the Appendix) for the number of components of $\mathbb{C} \setminus B[0, 1]$ of area greater than $\varepsilon > 0$ to the number of components of $\mathbb{C} \setminus B[0, n]$ of area greater than $n^{1-\delta}, \delta > 0$. Then couple planar random walk $S$ and planar Brownian motion $B$ via Skorokhod embedding and compare, for some $\delta_0 > 0$ and all $0 < \delta \leq \delta_0$, the number of holes of area larger than $n^{1-\delta}$ for $B[0, n]$ and $S[0, 2n]$ under the coupling, as described in the next steps.

2. Let $c = 1 + \varepsilon$, where $\varepsilon > 0$, and for $j \in \{0, 1, 2, \ldots\}, n \in \mathbb{N}$, define

(1.5) $$I_{j,n} = I_{j,n}(\delta, \varepsilon) = [n^{1-\delta}c^j, n^{1-\delta}c^{j+1}).$$

This gives a decomposition of $[n^{1-\delta}, \infty) = \bigcup_{j=0}^{\infty} I_{j,n}$. Let

$$N_{j,n} = N_{j,n}(\delta, \varepsilon) = \#\{\text{components of } \mathbb{C} \setminus B[0,n] \text{ with area in } I_{j,n}\},$$

$$\tilde{N}_{j,n} = \tilde{N}_{j,n}(\delta, \varepsilon) = \#\{\text{components of } \mathbb{C} \setminus S[0,2n] \text{ with area in } I_{j,n}\},$$

and show that for every $j \geq 0$, every $\varepsilon > 0$ small enough,

(1.6) $$\mathbb{P}(|\tilde{N}_{j,n} - N_{j,n}| > \varepsilon N_{j,n}) \to 0 \quad \text{fast enough, as } n \to \infty,$$

so that the sum over $j$ of these probabilities goes to 0, which then implies

$$\mathbb{P}(|\tilde{N}_n(\delta) - N_n(\delta)| > \varepsilon N_n(\delta)) \xrightarrow{n \to \infty} 0,$$

where

$$N_n(\delta) = \#\{\text{connected components of } \mathbb{C} \setminus B[0,n] \text{ of area } \geq n^{1-\delta}\}.$$

This is done by comparing the total area of all Brownian components with area in $I_{j,n}$ with the total area of all random walk components with area in $I_{j,n}$.

3. To show (1.6), analyze $\Delta(z) = \Delta_n(z) := ||C_n(z)| - |\tilde{C}_n(z)||$ for each $z \in \mathbb{Z}^2$, where $|C_n(z)|$ is the area of the connected component of $\mathbb{C} \setminus B[0,n]$ containing $z$, $|\tilde{C}_n(z)|$ is the area of the connected component of $\mathbb{C} \setminus S[0,2n]$ containing $z$. More specifically, show that on "good configurations," if $\max\{|C_n(z)|, |\tilde{C}_n(z)|\}$ is large, $\Delta(z)$ is of smaller order of magnitude than $\max\{|C_n(z)|, |\tilde{C}_n(z)|\}$.

4. Show that "bad configurations" are unlikely. This involves:

   - Handling the case where $z$ is close to $\partial C_n(z)$ or $\partial \tilde{C}_n(z)$, which is done with the help of ideas relating the two-sided disconnection exponent for Brownian motion and random walk to the fractal dimension of the Brownian frontier.
   - Looking at other "bad cases" which can occur even if $z$ is far from $\partial C_n(z)$ and $\partial \tilde{C}_n(z)$: $z$ being disconnected from $\infty$ by $B$ shortly before time $n$, not leaving $S$ enough time to disconnect $z$ from $\infty$ as well (or the other way around); $C_n(z)$ being a very thin component, and the same for $\tilde{C}_n(z)$. The tools used for these cases are the one-sided disconnection exponent and Beurling estimates for Brownian motion and random walk.

In Section 2, we give a list of definitions of the objects most commonly used throughout this paper and introduce Skorokhod embedding, the coupling which is at the center of our proof. In Section 3, we look at how "thick" the boundary of a Brownian motion or random walk component is. The ideas used are based on the method of [4], which exhibits the relationship between the two-sided disconnection exponent of Brownian motion and the Hausdorff dimension of the Brownian frontier. Section 4 contains a sequence of preparatory lemmas leading to a comparison between the areas of the Brownian motion and random walk component containing a given lattice point. In Section 5, we use the results of Sections 3 and 4 to prove Proposition 5.1 from which Theorem 1.1 follows immediately. In Section 6, we discuss the initial motivation of our study, namely Mandelbrot's observation of a different exponent for small lattice holes, and give a precise conjecture for this observation. In particular, we show that if the picture suggested by Mandelbrot is the right one, the regime in which the 5/3 exponent exists must be very small. Finally, the Appendix provides the consequences of [10] needed to make the link between small components of $\mathbb{C} \setminus B[0,1]$ and large components of $\mathbb{C} \setminus S[0,2n]$. Many basic estimates for planar Brownian motion and random walk are used throughout this paper, in particular large deviations and Beurling estimates. We will only include their statements in this paper, but their derivations, as well as a discussion on Skorokhod embedding, can be found in [1].

**2. Definitions and tools.** Throughout this paper, all multiplicative constants will be denoted by $K, K_1$ or $K_2$. It will be understood that they may be different from one line to the next. The letters $r, s, t$ will be used to denote real numbers, while $i, j, k, l, m, n$ will represent natural numbers. Points of the complex plane $\mathbb{C}$ will be represented by the letters $u, v, w, z$. The symbols $\mathcal{O}$, $o$ and $\sim$ will have the usual meaning: for two functions $f$ and $g$, $f(x) = \mathcal{O}(g(x))$ if there exists a constant $K$ such that $f(x) \leq Kg(x)$ for all $x$, $f(x) = o(g(x))$ if $\lim_{x \to \infty} f(x)/g(x) = 0$, and $f(x) \sim g(x)$ if $\lim_{x \to \infty} f(x)/g(x) = 1$. We say that a function $f(x)$ *decays rapidly* if for every $r \in \mathbb{R}, f(x) = o(x^{-r})$.

For any $z \in \mathbb{Z}^2$, we define $\text{Sq}(z)$ to be the closed region bounded by the square centered at $z$, whose sides are parallel to the axes and of length 1. For each $z \in \mathbb{Z}^2, \text{Sq}(z)$ will be called a *lattice square* (or just *square*). The Euclidean norm of a point $x$ in $\mathbb{R}$ or $\mathbb{C}$ is $|x|$ and the integer part of $x \in \mathbb{R}$ is $[x]$. The boundary of a set $A \subset \mathbb{C}$ will be denoted by $\partial A$, its area by $|A|$ and its diameter $\text{diam}(A) = \sup_{w,z \in A} |w - z|$. The distance between two sets $A, B \in \mathbb{C}$ is $d(A, B) = \inf_{x \in A, y \in B} |x - y|$. For any sets $A \subset F$, $A^c = F \setminus A$ will be the complement of $A$ in $F$. It will always be clear from context what $F$ is meant to be.

Unless stated otherwise, $B = (B(t))_{t \geq 0}$ will denote standard planar Brownian motion and $S = (S(n))_{n \geq 0}$ will stand for planar simple random walk. As a general rule, in this paper, a tilde will refer to a quantity related to random walk. $S$ will also denote the continuous process $(S(t))_{t \geq 0}$ obtained from planar simple random walk as in (1.1). It will be clear from the choice of the letter for the argument if we consider real or integer times.

At the center of the method used in this paper lies a coupling of planar random walk $S$ and planar Brownian motion $B$. It is a particular case of Skorokhod embedding. We state it here and refer the reader to [1] for more details.

PROPOSITION 2.1. *There exist a probability space $(\Omega, \mathcal{F}, \mathbb{P})$ containing a standard Brownian motion $B$ and simple random walk $S$ in the plane, constants $b, K > 0$, such that*

$$(2.1) \qquad \mathbb{P}\left(\sup_{0 \leq t \leq n} |B(t) - S(2t)| \geq n^{1/4} \log^2 n\right) \leq K n^{1 - b \log n}.$$

From here on, we will be working in this probability space. For any $t \geq 0$, we let

$$(2.2) \quad C_t(z) = \text{ the connected component of } \mathbb{C} \setminus B[0, t] \text{ containing } z$$

if $z \in \mathbb{C} \setminus B[0, t]$ and

$$\tilde{C}_t(z) = \text{ the connected component of } \mathbb{C} \setminus S[0, 2t] \text{ containing } z$$

if $z \in \mathbb{C} \setminus S[0, 2t]$. If $z \in B[0, t]$, we let $C_t(z) = \varnothing$ and if $z \in S[0, 2t], \tilde{C}_t(z) = \varnothing$. We will use the convention that $|\varnothing| = 0$. For any $z \in \mathbb{C}$ and $R \in \mathbb{R}_+$,

$$(2.3) \qquad D(z, R) = \{w \in \mathbb{C} : |w - z| \leq R\}$$

will denote the closed disk of radius $R$, centered at $z$ and $D(R)$ will be short for $D(0, R)$.

Two estimates for Brownian motion and random walk will be useful in the proofs of the next two propositions, as well as later in this paper. They are large deviations estimates giving an upper bound for the probability that in time $n$ random walk or Brownian motion travel much beyond distance $\sqrt{n}$ or remain in a disk of radius much smaller than $\sqrt{n}$. For the proofs, see [1].

LEMMA 2.1. *If $B$ is a planar Brownian motion, $S$ a planar simple random walk, there exists a constant $K > 0$ such that for every $n \geq 0$, every $r \geq 1$,*

$$(2.4) \qquad \mathbb{P}\left(\sup_{0 \leq t \leq n} |B(t)| \geq r\sqrt{n}\right) \leq K \exp\{-r^2/2\},$$

$$(2.5) \qquad \mathbb{P}\left(\max_{0 \leq k \leq 2n} |S(k)| \geq r\sqrt{n}\right) \leq K \exp\{-r^2/4\}.$$

LEMMA 2.2. *If $B$ is a planar Brownian motion, $S$ a planar simple random walk, there exists a constant $K > 0$ such that for every $n \geq 0$, every $r \geq 1$,*

(2.6)
$$\mathbb{P}\left(\sup_{0 \leq t \leq n} |B(t)| \leq r^{-1}\sqrt{n}\right) \leq \exp\{-Kr^2\},$$
$$\mathbb{P}\left(\max_{0 \leq k \leq 2n} |S(k)| \leq r^{-1}\sqrt{n}\right) \leq \exp\{-Kr^2\}.$$

We will also need on several occasions the following well known result which we give here without a proof (see [1] for the continuous case and [2] for the more difficult discrete case):

THEOREM 2.1 (Beurling estimate).

1. *There exists a constant $K > 0$ such that for any $R \geq 1$, any $x \in \mathbb{C}$ with $|x| \leq R$, any $A \subset \mathbb{C}$ with $[0, R] \subset \{|z| : z \in A\}$,*

(2.7)
$$\mathbb{P}^x(\xi_R \leq T_A) \leq K(|x|/R)^{1/2},$$

   *where $\xi_R = \inf\{t \geq 0 : |B(t)| \geq R\}$ and $T_A = \inf\{t \geq 0 : B(t) \in A\}$.*
2. *There exists a constant $K > 0$ such that for any $n \geq 1$, any $x \in \mathbb{Z}^2$ with $|x| \leq n$, any connected set $A \subset \mathbb{Z}^2$ containing the origin and such that $\sup\{|z| : z \in A\} \geq n$,*

$$\mathbb{P}^x(\Xi_n \leq \tau_A) \leq K(|x|/n)^{1/2},$$

   *where $\Xi_n = \inf\{k \geq 0 : |S(k)| \geq n\}$ and $\tau_A = \inf\{k \geq 0 : S(k) \in A\}$.*

**3. Disconnection exponents and the holes' boundaries.** Understanding the structure of the boundary of the components $C_n(z)$ and $\tilde{C}_n(z)$ containing a point $z \in \mathbb{C}$ is essential in our approach to the problem discussed in this paper. We first derive in Proposition 3.1 an upper bound, uniform for all $z \in \mathbb{Z}^2$, for the expected number of lattice squares (see Section 2) intersected by the boundary of $C_n(z)$ or $\tilde{C}_n(z)$. Proposition 3.2 then gives an upper bound for the expected number of lattice squares intersected by the union of the boundaries of all components of area greater than $n^{1-\delta}, \delta > 0$.

The main tool needed for our estimate is the two-sided disconnection exponent, one of many intersection and disconnection exponents computed exactly in [7, 8] and [9]. Since we will make use of the one-sided disconnection exponent in Section 4, we mention it here as well. It was shown in [5] and [6] that each of these exponents is the same for Brownian motion and random walk.

For any $x^1, x^2 \in \mathbb{C}$, we let $\mathbb{P}^{x^1, x^2}$ be the probability measure associated with two independent planar Brownian motions $B^1$ and $B^2$ with $B^1(0) = x^1$ and $B^2(0) = x^2$. $\mathbb{P}^{x^1}$ will denote the probability measure associated with the lone Brownian path $B^1$ started at $x^1$. For $i = 1, 2$ and $n \in \mathbb{N}$, we define $\Xi_n^i(x) = \inf\{t \geq 0 : |B^i(t) - x| \geq n\}$; if $x = 0$, we just write $\Xi_n^i$.

Since we will always deal separately with Brownian motion or random walk in this section, we can use the same notation for $S$ without risking any confusion: $\mathbb{P}^{x^1,x^2}$ is the probability measure associated to two independent planar random walks $S^1$ and $S^2$ with $S^1(0) = x^1$ and $S^2(0) = x^2$, where $x^1, x^2 \in \mathbb{Z}^2$ and $\mathbb{P}^{x^1}$ is the probability measure associated to $S^1$ started at $x^1$. It will be clear from context whether $\mathbb{P}$ refers to Brownian motion or random walk. We also let $\xi_n^i(x) = \inf\{k > 0 : |S^i(k) - x| \geq n\}$, and write $\xi_n^i = \xi_n^i(0)$, where $i = 1, 2$ and $n \in \mathbb{N}$. For any compact $A \subset \mathbb{C}$, we let

(3.1) $\quad \bar{Q}(A) = $ the closure of the unbounded component of $\mathbb{C} \setminus A$.

### 3.1. One-sided disconnection exponent. Let

$$A_n = \{D(1) \cap \bar{Q}(B^1[0, \Xi_n^1]) \neq \varnothing\} \quad \text{and} \quad \tilde{A}_n = \{(0,0) \cap \bar{Q}(S^1[0, \xi_n^1]) \neq \varnothing\},$$

where $D(1)$ is the closed unit disk centered at the origin, and write

$$\mathbb{P}(A_n) = \sup \mathbb{P}^x(A_n), \qquad \mathbb{P}(\tilde{A}_n) = \mathbb{P}^0(\tilde{A}_n),$$

where the sup is over all $x$ with $|x| \leq 1$.

The following lemma is a consequence of [7, 8] and [9], where the value of the Brownian disconnection exponent is computed and [5], where the equality between the Brownian and the random walk exponents is shown.

LEMMA 3.1. *There exists a constant $K > 0$ such that for all $n \geq 1$:*

(a) $\mathbb{P}(A_n) \leq K n^{-1/4}$,
(b) $\mathbb{P}(\tilde{A}_n) \leq K n^{-1/4}$.

### 3.2. Two-sided disconnection exponent. Let

$$F_n = \{D(1) \cap \bar{Q}(B^1[0,n] \cup B^2[0,n]) \neq \varnothing\},$$
$$D_n = \{D(1) \cap \bar{Q}(B^1[0, \Xi_n^1] \cup B^2[0, \Xi_n^2]) \neq \varnothing\},$$
$$\tilde{F}_n = \{(0,0) \cap \bar{Q}(S^1[0,n] \cup S^2[0,n]) \neq \varnothing\},$$
$$\tilde{D}_n = \{(0,0) \cap \bar{Q}(S^1[0, \xi_n^1] \cup S^2[0, \xi_n^2]) \neq \varnothing\}.$$

We will write

$$\mathbb{P}(F_n) = \sup \mathbb{P}^{x^1,x^2}(F_n) \quad \text{and} \quad \mathbb{P}(D_n) = \sup \mathbb{P}^{x^1,x^2}(D_n),$$

where the sup is over all $|x^1| \leq 1, |x^2| \leq 1$ and

$$\mathbb{P}(\tilde{F}_n) = \mathbb{P}^{0,0}(\tilde{F}_n) \quad \text{and} \quad \mathbb{P}(\tilde{D}_n) = \mathbb{P}^{0,0}(\tilde{D}_n).$$

Lemma 3.2 is based on the same papers as Lemma 3.1, except that the equality between the Brownian and random walk exponents follows from [6].

LEMMA 3.2. *There exists a constant $K > 0$ such that for all $n \geq 1$:*

(a) $\mathbb{P}(F_n) \leq K n^{-1/3}$,
(b) $\mathbb{P}(D_n) \leq K n^{-2/3}$,
(c) $\mathbb{P}(\tilde{F}_n) \leq K n^{-1/3}$,
(d) $\mathbb{P}(\tilde{D}_n) \leq K n^{-2/3}$.

3.3. *Two estimates for the holes' boundaries.* It was shown in [4] that there is a strong link between the two-sided disconnection exponent and the Hausdorff dimension of the Brownian frontier, defined as the boundary of the unbounded component of $\mathbb{C} \setminus B[0,1]$. This can be seen by observing that if $0 \leq t \leq 1$, $B(t)$ is in the frontier of $B[0,1]$ if $B[0,t] \cup B[t,1]$ does not disconnect $B(t)$ from infinity and that $B[t,1]$ and the time-reversal of $B[0,t]$, $(B(t-s))_{0 \leq s \leq t}$, are independent Brownian motions. [For two bounded sets $A, B \subset \mathbb{C}$, we will say that *$A$ does not disconnect $B$ from $\infty$ if $B \cap \bar{Q}(A) \neq \varnothing$*.] The proofs of the next two propositions are based on this idea.

In both propositions, the statements for random walk and Brownian motion are the same, but the fact that Brownian motion has more freedom to wander than random walk in a unit time interval makes the proofs slightly more technical in the Brownian case. For the sake of variety, we prove the first for random walk and the second for Brownian motion.

Proposition 3.1 gives an upper bound, uniform for all $z \in \mathbb{Z}^2$, for the expected number of lattice squares which are intersected by the boundary of $C_n(z)$ or $\tilde{C}_n(z)$. Recall that a function is rapidly decaying if it goes to 0 faster than any power function.

PROPOSITION 3.1. *There exists a constant $K > 0$ such that for every $z \in \mathbb{Z}^2$, every $n \geq 1$,*

$$\mathbb{E}\left[\sum_{y \in \mathbb{Z}^2} \mathbb{1}\{\mathrm{Sq}(y) \cap \partial C_n(z) \neq \varnothing\}\right] \leq K n^{2/3} (\log n)^{11/3}$$

*and*

(3.2) $$\mathbb{E}\left[\sum_{y \in \mathbb{Z}^2} \mathbb{1}\{\mathrm{Sq}(y) \cap \partial \tilde{C}_n(z) \neq \varnothing\}\right] \leq K n^{2/3} (\log n)^{11/3}$$

PROOF. We prove (3.2). By Lemma 2.1, there is a rapidly decaying function $\phi_1$ such that

$$\mathbb{E}\left[\sum_{y \in \mathbb{Z}^2} \mathbb{1}\{\mathrm{Sq}(y) \cap \partial \tilde{C}_n(z) \neq \varnothing\}\right]$$

$$\leq \sum_{j=0}^{2n} \mathbb{P}(S(j) \in \partial \tilde{C}_n(z))$$

(3.3)
$$= \sum_{w \in \mathbb{Z}^2} \sum_{j=0}^{2n} \mathbb{P}(S(j) \in \partial \tilde{C}_n(z); S(j) = w)$$

$$\leq \sum_{w \in \tilde{D}(\sqrt{n}\log n)} \sum_{j=0}^{2n} \mathbb{P}(S(j) \in \partial \tilde{C}_n(z); S(j) = w) + \phi_1(n),$$

where for $r \in \mathbb{R}_+, \tilde{D}(r) = D(r) \cap \mathbb{Z}^2 = \{z \in \mathbb{Z}^2 : |z| \leq r\}$. For all $0 \leq j \leq 2n$, we let

$$S_j^{(1)}(i) = S(j-i) - S(j), \qquad 0 \leq i \leq j,$$
$$S_j^{(2)}(i) = S(j+i) - S(j), \qquad 0 \leq i \leq 2n-j,$$

be the translates starting at the origin of the time-reversal of the part of $S$ up to time $j$ and the portion of the walk from time $j$ to $2n$, respectively. Then $S_j^{(1)}$ and $S_j^{(2)}$ are independent simple random walks, started at 0.

By (3.3), we just need to find a bound for every $z \in \mathbb{Z}^2$, $w \in \tilde{D}(\sqrt{n}\log n)$, and every $0 \leq j \leq 2n$, for

(3.4) $$\mathbb{P}(S(j) \in \partial \tilde{C}_n(z); S(j) = w).$$

To ease the reader's work, we quickly outline the ideas involved in the simplest case, where $|z| \geq \sqrt{n}\log n$ and $z$ is therefore very likely to be in $\bar{Q}(S[0, 2n])$, where $\bar{Q}$ is as in (3.1). In this case, for a typical $w$ and a typical walk $S$, the probability in (3.4) is bounded above by

$$\mathbb{P}(S(j) \in \bar{Q}(S[0,j] \cup S[j, 2n]); S(j) = w),$$

which is equal to

$$\mathbb{P}(0 \in \bar{Q}(S_j^{(1)}[0,j] \cup S_j^{(2)}[0, 2n-j]); S(j) = w).$$

Unfortunately, the events $\{0 \in \bar{Q}(S_j^{(1)}[0,j] \cup S_j^{(2)}[0, 2n-j])\}$ and $\{S(j) = w\}$ are not independent, since the first event influences the shape of $S_j^{(1)}[0,j]$ which has an influence on the second event. If they were, the local central limit theorem (see [1]) and Lemma 2.1 would imply that for $j \geq 1, \mathbb{P}(S(j) \in \partial \tilde{C}_n(z); S(j) = w) \leq K(\min\{j, n-j\})^{-1/3} \cdot j^{-1}$ and summing this over $1 \leq j \leq n-1$ and $w \in \tilde{D}(\sqrt{n}\log n)$ would give a bound of $Kn^{2/3}\log^2 n$. It turns out that this heuristic argument works if we consider the slightly different events

$$\{0 \in \bar{Q}(S_j^{(1)}[0, [j/\log^2 n]] \cup S_j^{(2)}[0, [(2n-j)/\log^2 n]])\}$$

and $\{S(j) = w\}$, which are "almost" independent. Indeed, $S_j^{(1)}[0, [j/\log^2 n]]$ is a much shorter path than $S_j^{(1)}[[j/\log^2 n], j]$ and its shape has little influence on the position of $S(j)$, which is mostly determined by $S_j^{(1)}[[j/\log^2 n], j]$.

The introduced logarithmic term causes almost no loss in the bound derived in the heuristic argument. For general $z$ and $w$, more work is required, especially since only the parts of $S_j^{(1)}$ and $S_j^{(2)}$ that lie in $D(w, |z-w|)$ determine whether $w \in \partial \tilde{C}_n(z)$ or not.

We need to estimate the probability in (3.4) in a different way in each of the following four exhaustive cases:

(i) $|z-w|^2 \leq j \leq n-j$,
(ii) $|z-w|^2 \leq n-j \leq j$,
(iii) $j \leq \min\{n-j, |z-w|^2\}$,
(iv) $n-j \leq \min\{j, |z-w|^2\}$.

We will consider cases (i) and (iii). Case (ii) [resp. (iv)] is handled exactly like that of (i) [resp. (iii)]. Note that $j \leq n-j \Leftrightarrow j \leq n/2$.

*Case* (i): For $w, z \in \mathbb{Z}^2$, we let $T = T(n, w, z) = [2(|z-w|/\log^2 n)^2]$ and define the event

$$\mathcal{T} = \mathcal{T}(j, n, w, z) = \left\{ \sup_{0 \leq i \leq T} |S_j^{(1)}(i)| \leq |z-w|; \sup_{0 \leq i \leq T} |S_j^{(2)}(i)| \leq |z-w| \right\}.$$

Note that by (2.5), $\phi_2(n) = \mathbb{P}(\mathcal{T}^c)$ decays rapidly and observe that if $S(j) = w$, $S(j) \in \partial \tilde{C}_n(z)$, and $\mathcal{T}$ occurs, then $S_j^{(1)}[0, T] \cup S_j^{(2)}[0, T]$ cannot disconnect 0 from infinity. We define

(3.5) $\qquad \mathcal{D}(A, B) = \{A \text{ does not disconnect } B \text{ from } \infty\},$

and will write $\mathcal{D} = \mathcal{D}(S_j^{(1)}[0, T] \cup S_j^{(2)}[0, T], 0)$. If $|z-w|^2 \leq j \leq n/2$, then

$$\mathbb{P}(S(j) \in \partial \tilde{C}_n(z); S(j) = w)$$
$$\leq \mathbb{P}(\mathcal{T}^c) + \mathbb{P}(\{S(j) = w\}; \mathcal{D}; \mathcal{T})$$
$$\leq \phi_2(n) + \sum \mathbb{P}(\{S(j-T) = u; S(j) = w\}; \mathcal{D}),$$

where the sum is over all $u \in \tilde{D}(w, |z-w|) = D(w, |z-w|) \cap \mathbb{Z}^2$ and $\phi_2(n)$ decays rapidly.

The Markov property applied at time $j - T$ and translation invariance of simple random walk can be used to see that this last expression is bounded above by

$$\sum \mathbb{P}(S(j-T) = u) \mathbb{P}(\mathcal{D}; \{S(T) = w - u\})$$
$$\leq \sup \mathbb{P}(S(j-T) = u) \sum \mathbb{P}(\mathcal{D}; \{S(T) = w - u\})$$
$$\leq \sup \mathbb{P}(S(j-T) = u) \mathbb{P}(\mathcal{D}),$$

where the sup and the sum are over the same set as above. Since $T \leq j/\log^4 n$, the local central limit theorem implies that the sup is bounded above by $K/j$ for some constant $K > 0$, independent of $j$ and $n$. Using Lemma 3.2(c) to bound $\mathbb{P}(\mathcal{D})$, we find that if $|z-w|^2 \leq j \leq n/2$, there exist a constant $K > 0$, independent of $z, w, j$ and $n$, and a rapidly decaying function $\phi_2$ such that

$$(3.6) \qquad \mathbb{P}(S(j) \in \partial \tilde{C}_n(z); S(j) = w) \leq \frac{K}{j}\left(\frac{|z-w|}{\log^2 n}\right)^{-2/3} + \phi_2(n).$$

*Case* (iii): If $j \leq \min\{n/2, |z-w|^2\}, U = U(j,n) = [\frac{j}{\log^4 n}]$, and

$$\mathcal{U} = \mathcal{U}(j,n,w,z) = \left\{\sup_{0 \leq i \leq U} |S_j^{(1)}(i)| \leq |z-w|; \sup_{0 \leq i \leq U} |S_j^{(2)}(i)| \leq |z-w|\right\},$$

then $\phi_3(n) = \mathbb{P}(\mathcal{U}^c)$ decays rapidly. Using the notation defined in (3.5), we write

$$\mathcal{D}' = \mathcal{D}(S_j^{(1)}[0,U] \cup S_j^{(2)}[0,U], 0).$$

Under $\mathcal{U}^c, \{S(j) \in \partial \tilde{C}_n(z)\}$ implies $\mathcal{D}'$. Therefore, using again the Markov property and the local central limit theorem,

$$(3.7) \quad \begin{aligned} &\mathbb{P}(S(j) \in \partial \tilde{C}_n(z); S(j) = w) \\ &\leq \sum \mathbb{P}(\{S(j) = w; S(j-U) = u\}; \mathcal{D}') + \phi_3(n) \\ &\leq \sup_{u \in \mathbb{Z}^2} \mathbb{P}(S(j-U) = u) \cdot \mathbb{P}(\mathcal{D}') + \phi_3(n) \\ &\leq \frac{K}{j}(j/\log^4 n)^{-1/3} + \phi_3(n) = K\left(\frac{j^2}{\log^2 n}\right)^{-2/3}. \end{aligned}$$

Here, the sum is over all $u \in \tilde{D}(w, \frac{|z-w|}{\log n})$ and $K$ is again independent of $z, w, j$, and $n$.

We find in the same way as in case (i) that if $|z-w|^2 \leq n - j \leq n/2$,

$$(3.8) \qquad \mathbb{P}(S(j) \in \partial \tilde{C}_n(z); S(j) = w) \leq \frac{K}{j}\left(\frac{|z-w|}{\log^2 n}\right)^{-2/3} + \phi_4(n)$$

and in the same way as in case (iii) that if $n - j \leq \min\{j, |z-w|^2\}$,

$$(3.9) \qquad \mathbb{P}(S(j) \in \partial \tilde{C}_n(z); S(j) = w) \leq K\left(\frac{(n-j)^2}{\log^2 n}\right)^{-2/3} + \phi_5(n).$$

We can now wrap up the proof by adding all our bounds: For a given $z$ with $|z-w|^2 \leq n/2$, we now obtain from (3.3), (3.6)–(3.9) that there exist a rapidly decaying function $\phi(n)$ and a constant $K$ such that for all $z \in \mathbb{Z}^2$ and $n \in \mathbb{N}$, with the notation $c = |z - w|^2$ and $\tilde{D} = \tilde{D}(\sqrt{n}\log n)$,

$$\mathbb{E}[\#\{y \in \mathbb{Z}^2 : \mathrm{Sq}(y) \cap \partial\tilde{C}_n(z) \neq \varnothing\}]$$

$$\leq \phi(n) + K \sum_{w \in \tilde{D}} \left[ 2 + \sum_{j=1}^{[c]} \left(\frac{j^2}{\log^2 n}\right)^{-2/3} \right.$$

$$\left. + \sum_{j=[c \wedge n/2]}^{n-[c \wedge n/2]} \frac{1}{j}\left(\frac{|z-w|}{\log^2 n}\right)^{-2/3} + \sum_{j=n-[c]}^{n-1} \left(\frac{(n-j)^2}{\log^2 n}\right)^{-2/3} \right]$$

$$\leq \phi(n) + K \log^{4/3} n \sum_{w \in \tilde{D}(\sqrt{n}\log n)} \left( 2 + 2\sum_{j=1}^{[c]} j^{-4/3} + \sum_{j=[c \wedge n/2]}^{[n/2]} \frac{|z-w|^{-2/3}}{j} \right)$$

$$\leq \phi(n) + K \log^{4/3} n \sum_{w \in \tilde{D}(\sqrt{n}\log n)} |z-w|^{-2/3}(1 + \log n)$$

$$\leq K n^{2/3} (\log n)^{11/3}. \qquad \square$$

For $0 \leq \delta \leq 1$, we define

$$\mathcal{C}_n = \mathcal{C}_n(\delta) = \{z \in \mathbb{C} : |C_n(z)| \geq n^{1-\delta}\},$$
$$\tilde{\mathcal{C}}_n = \tilde{\mathcal{C}}_n(\delta) = \{z \in \mathbb{C} : |\tilde{C}_n(z)| \geq n^{1-\delta}\}.$$

In the proof of Proposition 5.1, we will need an estimate for the expected number of lattice squares (for the definition, see Section 2) which are intersected by the boundary of $\mathcal{C}_n(\delta)$ and $\tilde{\mathcal{C}}_n(\delta)$. The ideas of the proof below are similar to those of Proposition 3.1, but the fact that we do not have to worry about the location of $S(j)$ relatively to a point $z$ makes this next proof considerably simpler.

PROPOSITION 3.2. *There exists a constant $K > 0$ such that for every $\delta > 0$, every $n \geq 1$,*

(3.10) $\quad \mathbb{E}[\#\{y \in \mathbb{Z}^2 : \mathrm{Sq}(y) \cap \partial\mathcal{C}_n(\delta) \neq \varnothing\}] \leq K n^{(2+\delta)/3} (\log n)^{8/3}$

*and*

(3.11) $\quad \mathbb{E}[\#\{y \in \mathbb{Z}^2 : \mathrm{Sq}(y) \cap \partial\tilde{\mathcal{C}}_n(\delta) \neq \varnothing\}] \leq K n^{(2+\delta)/3} (\log n)^{8/3}.$

PROOF. We prove (3.10) which presents a minor additional difficulty, since we have less control over the diameter of a path of time length one for Brownian motion than for random walk. The other ideas for the proof of (3.11) are the same as for (3.10). The general strategy of the proof is to find bounds for the expected number of time segments $[j-1, j]$ over which the Brownian path intersects $\partial\mathcal{C}_n(\delta)$. This will suffice since the expected number of lattice squares intersected by $B[j-1,j]$ is finite. For $1 \leq j \leq n$, let $d_j = \mathrm{diam}(B[j-1,j])$ and $\mathcal{B}_j = D(B(j), d_j)$ be the closed disk of radius $d_j$, centered at $B(j)$, so that $B[j-1,j] \subset \mathcal{B}_j$. We let $\bar{n} = [n^{1-\delta} \log^{-4} n]$ and define

$$\mathcal{A}_j^L = B[\max\{0, j-1-\bar{n}\}, j-1] \quad \text{and} \quad \mathcal{A}_j^R = B[j, \min\{2j-1, j+\bar{n}\}].$$

$\mathcal{A}_j^L$ and $\mathcal{A}_j^R$ span time intervals of same length and that length is $j-1$ if $j \leq \bar{n}+1$ and $\bar{n}$ if $j \geq \bar{n}+1$.

Our choice of $\bar{n}$ is motivated by the following: By Lemma 2.1, it is very likely that for every $1 \leq j \leq n$, $\mathcal{A}_j^L \cup \mathcal{A}_j^R$ is completely contained in a disk of radius $n^{(1-\delta)/2}/2$. More precisely, if $\mathcal{V} = \{\mathcal{A}_j^L \cup \mathcal{A}_j^R \subset D(B(j), n^{(1-\delta)/2}/2)\}$, then $\mathbb{P}(\mathcal{V}^c) = \phi(n)$, where $\phi$ is a rapidly decaying function. On the event $\mathcal{V}$, then, if $B[j-1,j]$ intersects $\partial\mathcal{C}_n(\delta)$, $\mathcal{A}_j^L \cup \mathcal{A}_j^R$ cannot disconnect $\mathcal{B}_j$ from infinity. Indeed, if it did, $\mathcal{B}_j$ would be disconnected from infinity by a portion of the Brownian path entirely contained in a disk of radius $n^{(1-\delta)/2}/2$ and could intersect the boundary of components of area no more than $n^{1-\delta}$.

The rest of the argument just involves taking care of the case where $d_j$ is unusually large and applying Lemma 3.2(a). Let $m = [n/2] + 1$. By symmetry, we have

$$(3.12) \quad \begin{aligned} &\mathbb{E}[\#\{y \in \mathbb{Z}^2 : \mathrm{Sq}(y) \cap \partial\mathcal{C}_n \neq \varnothing\}] \\ &= 2 \sum_{j=1}^m \mathbb{E}[\#\{y \in \mathbb{Z}^2 : \mathrm{Sq}(y) \cap \partial\mathcal{C}_n \cap B[j-1,j] \neq \varnothing\}]. \end{aligned}$$

We first look at the terms for $1 \leq j \leq \bar{n}$ (in which case $\mathcal{A}_j^L$ and $\mathcal{A}_j^R$ each span a time interval of length $j$). There is a $K > 0$ such that

$$\mathbb{E}[\#\{y \in \mathbb{Z}^2 : \mathrm{Sq}(y) \cap \partial\mathcal{C}_n \cap B[j-1,j] \neq \varnothing\}]$$

$$= \sum_{l \geq 1} \mathbb{E}[\#\{y \in \mathbb{Z}^2 : \mathrm{Sq}(y) \cap \partial\mathcal{C}_n \cap B[j-1,j] \neq \varnothing\}; d_j \in [l-1, l]]$$

$$\leq K \sum_{l \geq 1} l^2 \mathbb{P}(\partial\mathcal{C}_n \cap \mathcal{B}_j \neq \varnothing; d_j \in [l-1, l])$$

$$(3.13) \quad \leq K \left( \log^2 n (\mathbb{P}(\{\mathcal{B}_j \in \bar{Q}(\mathcal{A}_j^L \cup \mathcal{A}_j^R); d_j \leq \log n\}; \mathcal{V}) + \phi(n)) \right.$$

$$\left. + \sum_{l > \log n} l^2 \mathbb{P}(d_j \in [l-1, l]) \right)$$

$$\leq K \left( \log^2 n \left( \frac{j}{\log^2 n} \right)^{-1/3} + \phi(n) + \sum_{l > \log n} l^2 \exp\left\{ -\frac{(l-1)^2}{2} \right\} \right)$$

$$\leq K j^{-1/3} \log^{8/3} n,$$

where $\phi$ decays rapidly and the last inequality follows from Brownian scaling and Lemma 2.1. Therefore,

$$\sum_{j=1}^{\bar{n}} \mathbb{E}[\#\{y \in \mathbb{Z}^2 : \mathrm{Sq}(y) \cap \partial \mathcal{C}_n \cap B[j-1, j] \neq \varnothing\}]$$

$$\leq K(\log n)^{8/3} \sum_{j=1}^{\bar{n}} j^{-1/3}$$

$$\leq K n^{(2-2\delta)/3}.$$

In the same way, we find

$$\sum_{j=\bar{n}+1}^{m} \mathbb{E}[\#\{y \in \mathbb{Z}^2 : \mathrm{Sq}(y) \cap \partial \mathcal{C}_n \cap B[j-1, j] \neq \varnothing\}]$$

$$\leq K n \log^2 n \left(\frac{\log^2 n}{\bar{n}}\right)^{1/3}$$

$$\leq K n^{(2+\delta)/3} \log^{8/3} n,$$

which, combined with (3.12) and (3.13), gives the proposition. □

**4. Comparing Brownian and random walk areas.** A key step in the proof of Theorem 1.1 is to show that if $B$ and $S$ are coupled as in Proposition 2.1, then under certain favorable conditions, the difference between the areas of the Brownian and the random walk components containing a given point $z \in \mathbb{Z}^2$, $\Delta_n(z) = ||C_n(z)| - |\tilde{C}_n(z)||$, has an expected value of smaller order of magnitude than the areas themselves. The conditions are that at least one of $|C_n(z)|$ and $|\tilde{C}_n(z)|$ is $\geq n^{1-b}$ for some specific $b \in (0, 1)$, and that $z$ is not too close to $\partial C_n(z)$. This estimate is given in Proposition 4.1 and requires finding an upper bound, for fixed $z$, for the number of points $y \in \mathbb{Z}^2$ with $y \in (C_n(z) \setminus \tilde{C}_n(z)) \cup (\tilde{C}_n(z) \setminus C_n(z))$. Of course this number is only meaningful if $|C_n(z)| < \infty$ and $|\tilde{C}_n(z)| < \infty$ and we will estimate $\mathbb{E}[\Delta_n(z)]$ with the assumption that this is the case. The points for which $|C_n(z)| = \infty$ or $|\tilde{C}_n(z)| = \infty$ will be dealt with in the proof of Proposition 5.1. A useful estimate for that purpose is derived in Lemma 4.3.

If both components are finite, a point $y$ is in $C_n(z) \setminus \tilde{C}_n(z)$ if it has been disconnected from $z$ by $S$ but not by $B$. Informally, this can happen if:

1. Either $y$ or $z$ is closer to $\partial C_n(z)$ than the distance prescribed by the coupling in (2.1).

2. $y$ is disconnected from $z$ by $S$ very late, that is, at a time close to $2n$, in which case, $B$ gets "very close" to disconnecting $y$ from $z$, but may not have time to do so.
3. $B$ gets very close to disconnecting $y$ from $z$, but does not, despite having plenty of time to do so.

The main results needed to handle Case 1. were derived in the previous section. We provide the estimates needed for Cases 2 and 3 in Lemmas 4.5 and 4.4, respectively. These two lemmas are then used to prove Lemma 4.6, which gives a bound for the probability that $y \in (C_n(z) \setminus \tilde{C}_n(z)) \cup (\tilde{C}_n(z) \setminus C_n(z))$ if both components are bounded, one of them is large enough, and $y, z$ are not too close to $\partial C_n(z)$ or $\partial \tilde{C}_n(z)$.

Lemmas 4.1 and 4.2 give estimates similar to those in Lemmas 4.4 and 4.5, but are concerned with the first time a point $z$ is disconnected from $\infty$ rather than from another point $y$. They are used to prove Lemma 4.3, which, as mentioned above, is used to deal with the case where $|C_n(z)| = \infty$ or $|\tilde{C}_n(z)| = \infty$.

Note that Lemmas 4.1 and 4.4 are results about Brownian motion or random walk only, while the other lemmas of this section address questions about the joint behavior of the coupled random walk and Brownian motion. All the lemmas in this section have statements that come in pairs, where one version is obtained from the other by interchanging the roles of $B$ and $S$.

We will be interested in the first time at which an arbitrary point lies in a finite component and in studying how long it takes from a time at which it "almost" lies in a finite component until it actually does. This motivates the definition of the *closing times* for $z$ by Brownian motion and random walk:

$$T_z = \inf\{t \geq 0 : 0 < |C_t(z)| < \infty\}, \qquad \tilde{T}_z = \inf\{t \geq 0 : 0 < |\tilde{C}_t(z)| < \infty\},$$

where $C_t(z)$ [resp. $\tilde{C}_t(z)$] is the connected component of $\mathbb{C} \setminus B[0,t]$ [resp. $\mathbb{C} \setminus S[0,2t]$] containing $z$ and $|C_t(z)|$ [resp. $|\tilde{C}_t(z)|$] was defined to be 0 if $z \in B[0,t]$ (resp. $z \in S[0,2t]$). We use here the convention that $\inf \varnothing = \infty$. The *closing points* for $z$ are

$$x_z = B(T_z) \quad \text{and} \quad \tilde{x}_z = S(\tilde{T}_z).$$

For a curve $\gamma : [a,b] \to \mathbb{C}$ and a point $z \notin \{\gamma(t)\}_{a \leq t \leq b}$, we let $\arg_z(\gamma(t))$ denote the continuous argument of $\gamma$ about $z$, with the convention that $\arg_z(\gamma(0)) = 0$. This is well defined (see [13]). Note that $\arg_z(\cdot)$ is defined on the parametric interval $[a,b]$, not on the image of $\gamma$. The curves $\gamma$ with which we will be working are of course $B$ and $S$ and as the points $z$ we are interested in will always lie off the paths of $B$ and $S$, the argument will always be well defined. We assume henceforth that this is the case. We will use the abbreviation

$$\arg_z(\gamma(t), \gamma(s)) = \arg_z(\gamma(t)) - \arg_z(\gamma(s)).$$

This definition allows us to give another characterization of $T_z$:

$$T_z = \inf\{t \geq 0 : \exists 0 \leq s < t \text{ with } B(s) = B(t), |\arg_z(B(s), B(t))| \neq 0\},$$

and similarly for $\tilde{T}_z$. Note that $T_z$ and $\tilde{T}_z$ are stopping times. The *last call* for $z$ by B and the *last call* for $z$ by S are, respectively,

$$T_z^l = T_z^l(n) = \inf\left\{t \geq 0 : \exists s \in [0, t] : |\arg_z(B(s), B(t))| \geq \frac{3\pi}{2};\right.$$

$$\left. |B(s) - B(t)| \leq 3n^{1/4} \log^2 n\right\},$$

$$\tilde{T}_z^l = \tilde{T}_z^l(n) = \inf\left\{t \geq 0 : \exists r \in [0, t] : |\arg_z(S(2r), S(2t))| \geq \frac{3\pi}{2};\right.$$

$$\left. |S(2r) - S(2t)| \leq 3n^{1/4} \log^2 n\right\}.$$

Observe the factors of 2 in the definition of $\tilde{T}_z^l$, which are due to the fact that in the Skorokhod embedding of Proposition 2.1, $B$ and $S$ run on different clocks. The *last call points* are $x_z^l = B(T_z^l)$ and $\tilde{x}_z^l = S(2\tilde{T}_z^l)$. Note that by continuity of $B$ and $S$, $T_z^l < T_z$ and $\tilde{T}_z^l < \tilde{T}_z$.

The first lemma of this section shows that typically the closing time by $B$ for any point comes "soon" after the last call. For such a point $z$, and any $a > 1/2$, this lemma will imply that if $T_z^l \leq n - n^a$, it is unlikely that $z$ will not be disconnected from infinity by $B[0, n]$. Since if $z$ is not too close to $\partial C_n(z)$ and $B$ and $S$ are close to each other, we have $T_z^l \leq \tilde{T}_z$, this means that $\{\tilde{C}_n(z) < \infty; C_n(z) = \infty\}$ is an unlikely event unless $\tilde{T}_z \geq n - n^a$, that is, $z$ is disconnected from $\infty$ by $S$ "very late." This in turn will be shown to be unlikely in Lemma 4.2.

LEMMA 4.1. *There exists a constant $K > 0$ such that for any $z \in \mathbb{Z}^2$, any $n \in \mathbb{N}$ and any $a > 1/2$:*

(a) $\mathbb{P}(T_z - T_z^l > n^a) \leq Kn^{(1-2a)/24} \log n$.
(b) $\mathbb{P}(\tilde{T}_z - \tilde{T}_z^l > n^a) \leq Kn^{(1-2a)/24} \log n$.

PROOF. (a) To find a bound for $\mathbb{P}(T_z - T_z^l > n^a)$ we will need to consider all the relative positions of $z$ and $x_z^l = B(T_z^l)$. This will yield two different bounds. The first will be better for larger $|z - x_z^l|$, while the second will be better for smaller $|z - x_z^l|$.

First suppose that $m(z,a) := \min\{|z - x_z^l|, n^{a/2}/\log n\} \geq 300n^{1/4}\log^2 n$. Note that if $u \in D(z, m(z,a)/2)^c$ [recall the definition of $D$ in (2.3)] and $0 \leq s \leq t$ are such that $B[s,t] \subset D(u, m(z,a)/100)$, then, for all $n \geq 2$, we have the obvious rough bound

$$(4.1) \qquad |\arg_z(B(s), B(t))| \leq \frac{\pi}{8}.$$

We define $T_z^f = \inf\{t \geq 0 : |\arg_z(B(t), B(T_z^l))| \geq 3\pi/2, |B(t) - x_z^l| \leq 3n^{1/4} \times \log^2 n\}$, $x_z^f = B(T_z^f)$, and $\Phi_z^f = \inf\{t \geq T_z^f : B(t) \in \partial D_z^f\}$, where $D_z^f = D(x_z^f, m(z,a)/100)$. The point $x_z^f$ can be thought of as lying "across from" $x_z^l$ on $B[0, T_z^l]$ (see Figure 2). We first note that the definition of $x_z^l$ implies that the connected random set $A_z = B[T_z^f, \Phi_z^f]$ contains $x_z^f$, intersects $\partial D_z^f$, and satisfies:

1. $d(x_z^l, A_z) \leq 3n^{1/4}\log^2(n)$.
2. For any $t \in [T_z^f, \Phi_z^f], |\arg_z(B(t), B(T_z^l))| > \pi$.

The second point is true because $|\arg_z(B(T_z^f), B(T_z^l))| \geq 3\pi/2$ and inside $D_z^f, \arg_z(B(t))$ does not vary by more than $\frac{\pi}{8}$, by (4.1), since $x_z^f \in D(z, m(z,a)/2)^c$.

If we let $\Phi_z^l = \inf\{t \geq T_z^l : B(t) \in \partial D_z^f\}$, then

$$\mathbb{P}(T_z - T_z^l > n^a) \leq \mathbb{P}(B[T_z^l, T_z^l + n^a] \cap A_z = \varnothing)$$
$$\leq \mathbb{P}(B[T_z^l, T_z^l + n^a] \cap \partial D_z^f = \varnothing)$$
$$+ \mathbb{P}(B[T_z^l, \Phi_z^l] \cap A_z = \varnothing).$$

Therefore, using the fact that $T_z^l$ is a stopping time, we can combine the fact that

$$\phi(n) = \mathbb{P}(B[T_z^l, T_z^l + n^a] \cap \partial D_z^f = \varnothing) \leq \mathbb{P}\left(\sup_{0 \leq t \leq n^a} |B(t) - B(0)| \leq m(z,a)\right)$$

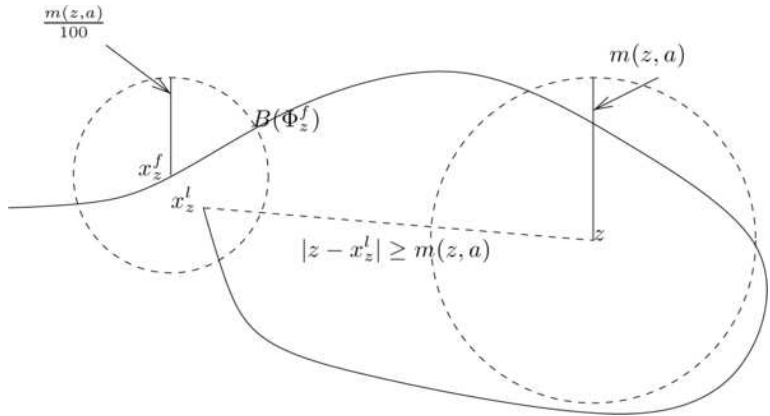

FIG. 2. *Lemma 4.1: $T_z^l$ is the first time at which $z$ is "almost" disconnected from $\infty$.*

decays rapidly (2.6) and the Beurling estimate (2.7) to obtain

$$\mathbb{P}(T_z - T_z^l > n^a) \leq \phi(n) + K\left(\frac{n^{1/4}\log^2 n}{m(z,a)}\right)^{1/2}$$
(4.2)
$$\leq K n^{1/8} m(z,a)^{-1/2} \log n.$$

We now derive a second bound in the case where $m(z,a) = |z - x_z^l|$. If we let $\hat{D} = D(x_z^l, n^{a/2}/\log n)$, define $\hat{\Xi} = \inf\{t \geq T_z^l : B(t) \in \hat{D}\}$, and use the definition of $\bar{Q}$ in (3.1), we find that since $z \in D(x_z^l, |z - x_z^l|)$, there is a rapidly decaying function $\phi$ such that

$$\begin{aligned}
\mathbb{P}(T_z - T_z^l > n^a) &\leq \mathbb{P}(z \in \bar{Q}(B[T_z^l, T_z^l + n^a])) \\
&\leq \mathbb{P}(D(x_z^l, |z - x_z^l|) \cap \bar{Q}(B[T_z^l, T_z^l + n^a]) \neq \varnothing) \\
&\leq \mathbb{P}(D(x_z^l, |z - x_z^l|) \cap \bar{Q}(B[T_z^l, \hat{\Xi}]) \neq \varnothing) + \mathbb{P}(\hat{\Xi} - T_z^l > n^a) \\
&\leq K\left(\frac{|z - x_z^l|}{n^{a/2}/\log n}\right)^{1/4} + \phi(n) \\
&\leq K n^{-a/8} |z - x_z^l|^{1/4} \log^{1/4} n,
\end{aligned}$$
(4.3)

by (2.6) and a scaled version of Lemma 3.1(a).

We can now conclude by noting that if $m(z,a) = |z - x_z^l| \geq 300 n^{1/4} \log^2 n$, we have two bounds to choose from and find from (4.2) and (4.3) that

(4.4) $\quad \mathbb{P}(T_z - T_z^l > n^a; m(z,a) = |z - x_z^l|) \leq K n^{(1-2a)/24} \log n.$

Also, if $m(z,a) = n^{a/2}/\log n$, (4.2) yields

(4.5) $\quad \mathbb{P}(T_z - T_z^l > n^a; m(z,a) = n^{a/2}/\log n) \leq K n^{(1-2a)/8} \log^{3/2} n.$

Combining (4.4) and (4.5) gives the lemma. □

There are obviously many times in the interval $[0,n]$ at which new points become disconnected from infinity by $B$. However, for any given point $z \in \mathbb{Z}^2$, the probability that $z$ becomes disconnected from $\infty$ "late" in the interval $[0,n]$ is small.

LEMMA 4.2. *There exists a constant $K > 0$ such that for any $z \in \mathbb{Z}^2, n \in \mathbb{N}$ and $a > 1/2$:*

(a) $\mathbb{P}(T_z \in [n - n^a, n]) \leq K n^{(a-1)/3} \log^{10/3} n,$
(b) $\mathbb{P}(T_z^l \in [n - n^a, n]) \leq K n^{(a-1)/3} \log^{10/3} n,$
(c) $\mathbb{P}(\tilde{T}_z \in [n - n^a, n]) \leq K n^{(a-1)/3} \log^{10/3} n,$
(d) $\mathbb{P}(\tilde{T}_z^l \in [n - n^a, n]) \leq K n^{(a-1)/3} \log^{10/3} n.$

PROOF. We prove (b). The other cases use exactly the same ideas. For $w \in \mathbb{C}$, let $\hat{D}_w = \hat{D}_w(a,n) := D(w, 4n^{a/2} \log^2 n)$ and recall the definitions made in the proof of the previous lemma: $T_z^f = \inf\{t \geq 0 : |\arg_z(B(t), B(T_z^l))| \geq 3\pi/2, |B(t) - x_z^l| \leq 3n^{1/4} \log^2 n\}$, $x_z^f = B(T_z^f)$. Since $a > 1/2$, if $x_z^l \in D(w, n^{a/2})$, then $x_z^f \in \hat{D}_w$. Therefore, using Lemma 2.1,

$$\mathbb{P}(T_z^l \in [n - n^a, n])$$
(4.6)
$$\leq \sum \mathbb{P}(T_z^l \in [n - n^a, n]; x_z^l \in D(w, n^{a/2})) + \phi(n)$$
$$= \sum \mathbb{P}(T_z^l \in [n - n^a, n]; x_z^l \in D(w, n^{a/2}); x_z^f \in \hat{D}_w) + \phi(n),$$

where the sum is over $w \in (n^{a/2} \cdot \mathbb{Z}^2) \cap D(\sqrt{n} \log n)$, $\phi$ decays rapidly and

$$n^{a/2} \cdot \mathbb{Z}^2 = \{(x,y) \in \mathbb{R}^2 : x = kn^{a/2}, y = ln^{a/2}, k, l \in \mathbb{Z}\}.$$

The definition of $T_z^l$ implies that if $a > 1/2$, then

(4.7) $\{T_z^l \in [n - n^a, n]; x_z^l \in D(w, n^{a/2}); x_z^f \in \hat{D}_w\} \subset \{z \notin \bar{Q}(B[0, T_z^l] \cup \hat{D}_w)\},$

where $\bar{Q}$ is defined as in (3.1). (See Figure 3.) Define $\Phi = \Phi_w = \inf\{t \geq 0 : B(t) \in \hat{D}_w\}$, $\Gamma = \Gamma_{z,w} = \inf\{t \geq \Phi : |B(t) - w| = |z - w|\}$, and $\Gamma^l = \Gamma^l_{z,w} = \sup\{t \leq T_z^l : |B(t) - w| = |z - w|\}$. The event

$$\{z \notin \bar{Q}(B[0, T_z^l] \cup \hat{D}_w)\} \cap \{\bar{Q}(B[\Phi, \Gamma] \cup B[\Gamma^l, T_z^l]) \cap \hat{D}_w = \varnothing\}$$

is contained in

$$\{z \notin \bar{Q}(B[0, T_z^l] \cup \hat{D}_w)\} \cap \{\bar{Q}(B[0, T_z^l]) \cap \hat{D}_w = \varnothing\} \subset \{z \notin \bar{Q}(B[0, T_z^l])\}.$$

But the continuity of the Brownian path implies that $T_z^l < T_z$, so this last event is the empty set. It follows from (4.7) that with the definition $A_{z,w} = A_{z,w}(n) = \{x_z^l \in D(w, n^{a/2}); x_z^f \in \hat{D}_w\}$,

$$\{T_z^l \in [n - n^a, n]; A_{z,w}\} \subset \{\bar{Q}(B[\Phi, \Gamma] \cup B[\Gamma^l, T_z^l]) \cap \hat{D}_w \neq \varnothing; A_{z,w}\}.$$

Therefore, (4.6) becomes

$$\mathbb{P}(T_z^l \in [n - n^a, n])$$
$$\leq \sum \mathbb{P}(\{\bar{Q}(B[\Phi, \Gamma] \cup B[\Gamma^l, T_z^l]) \cap \hat{D}_w \neq \varnothing\}; A_{z,w}) + \phi(n).$$

The two events in this last probability are not independent, but, exactly as in Proposition 3.1, we can write the probability of their intersection as the product of their probabilities, up to a logarithmic correction term. Noting that $B[\Phi, \Gamma]$ and the time-reversal of $B[\Gamma^l, T_z^l]$ are two independent Brownian paths from the inside of $\hat{D}_w$ to $D(w, |z - w|)$ and that by the local central limit theorem (see [3]), $\mathbb{P}(x_z^l \in D(w, n^{a/2}); T_z^l \geq n - n^a) \leq K n^{a-1}$, where $K$ is uniform for all $a > 1/2$, we can use Lemma 3.2 to find that

$$\mathbb{P}(T_z^l \in [n - n^a, n]) \leq K \sum \left( \frac{|z-w|}{n^{a/2} \log^3 n} \right)^{-2/3} n^{a-1},$$

where the sum is over $w \in (n^{a/2} \cdot \mathbb{Z}^2) \cap D(\sqrt{n} \log n)$. This can easily be seen to be bounded above by $K n^{(a-1)/3} \log^{10/3} n$. □

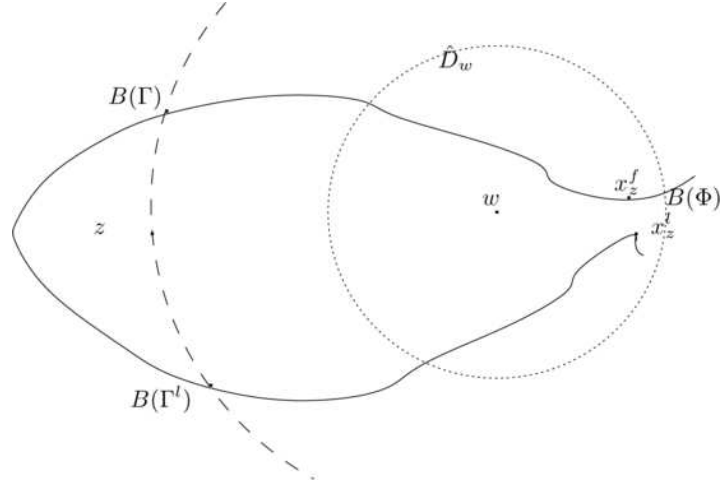

FIG. 3. *Lemma* 4.2: $B[0, T_z^l] \cup \hat{D}$ *disconnects* $z$ *from infinity. The portion of the Brownian path inside* $D(w, |z-w|)$ *cannot disconnect* $\hat{D}$ *from infinity.*

Lemma 4.3 shows that if $B$ and $S$ are coupled as in Proposition 2.1 the chance that $C_n(z)$ is finite while $\tilde{C}_n(z)$ is infinite is small when $z$ is not too close to the boundary of $C_n(z)$. The condition for $z$ to be away from the boundary of $C_n(z)$ is essential, since otherwise, we have no control over the probability that the Brownian path passes on one "side" of $z$ and the random walk on the other. To avoid this, we will choose $z$ to be at a distance from $\partial C_n(z)$ which is greater than the maximal distance between the coupled random walk and Brownian motion, that is, $n^{1/4} \log^2 n$. We define for $z \in \mathbb{C}$ the events

$$
\begin{aligned}
\mathcal{B}(z) &= \{d(z, \partial C_n(z)) \geq 100 n^{1/4} \log^2 n\}, \\
\tilde{\mathcal{B}}(z) &= \{d(z, \partial \tilde{C}_n(z)) \geq 100 n^{1/4} \log^2 n\},
\end{aligned}
$$
(4.8)

for which the dependence on $n$ should be noted but will not be written explicitly in order to keep the notation simple. As always $d(\cdot, \cdot)$ denotes Euclidean distance. We also define

$$\mathcal{P} = \mathcal{P}_n = \left\{ \sup_{0 \leq t \leq n} |B(t) - S(2t)| \leq n^{1/4} \log^2 n \right\},$$

the condition that $B$ and $S$ are close to each other as in the Skorokhod embedding. By (2.1) $\mathbb{P}(\mathcal{P}_n^c)$ decays rapidly.

LEMMA 4.3. *There exists a constant $K > 0$ such that for every $z \in \mathbb{Z}^2$ and $n \in \mathbb{N}$:*

(a) $\mathbb{P}(\{|C_n(z)| < \infty; |\tilde{C}_n(z)| = \infty\}; \mathcal{B}(z); \mathcal{P}) \leq Kn^{-1/30} \log^{10/3} n$.
(b) $\mathbb{P}(\{|\tilde{C}_n(z)| < \infty; |C_n(z)| = \infty\}; \tilde{\mathcal{B}}(z); \mathcal{P}) \leq Kn^{-1/30} \log^{10/3} n$.

PROOF (See Figure 4). (a) The idea of the proof is that on the event $\mathcal{B}(z)$, at the instant $z$ is disconnected from infinity by $B$, $\tilde{T}_z^l$, the last call for $z$ by $S$ has already occurred. At that instant, either the time is very close to $n$, which is unlikely by Lemma 4.2(a), or the random walk has plenty of time to disconnect $z$ from infinity and will do it with high probability by Lemma 4.1(b).

We use the fact that for $n$ large enough, $\mathcal{B}(z) \cap \mathcal{P} \subset \{\tilde{T}_z^l \leq T_z\}$, which follows from an argument similar to the one we used at the beginning of the proof of Lemma 4.1: By definition of $T_z$, there is a time $T \leq T_z$ such that $B(T) = B(T_z)$ and $|\arg_z(B(T), B(T_z))| \geq 2\pi$. The conditions $\mathcal{B}(z)$ and $\mathcal{P}$ ensure that $|\arg_z(B(T)) - \arg_z(S(2T))| \leq \frac{\pi}{8}$ and $|\arg_z(B(T_z)) - \arg_z(S(2T_z))| \leq \frac{\pi}{8}$, so that

$$(4.9) \qquad |\arg_z(S(2T), S(2T_z))| \geq \frac{3\pi}{2}.$$

Also, the condition $\mathcal{P}$ guarantees that $|S(2T) - B(T)| \leq n^{1/4} \log^2 n$ and $|S(2T_z) - B(T_z)| \leq n^{1/4} \log^2 n$, which implies

$$(4.10) \qquad |S(2T) - S(2T_z)| \leq 2n^{1/4} \log^2 n.$$

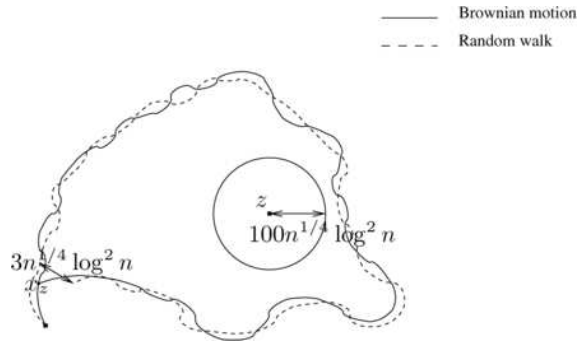

FIG. 4. *Proof of Lemma 4.3: when Brownian motion disconnects $z$ from $\infty$, random walk is close to doing so as well.*

From the definition of $\tilde{T}_z^l$, (4.9) and (4.10), it is now clear that $\tilde{T}_z^l \leq T_z$. Therefore, we have, for any $z \in \mathbb{C}$ and $n$ large enough,

$$\mathbb{P}(\{|C_n(z)| < \infty; |\tilde{C}_n(z)| = \infty\}; \mathcal{B}(z); \mathcal{P})$$
$$\leq \mathbb{P}(|\tilde{C}_n(z)| = \infty; \tilde{T}_z^l < n - n^{9/10}) + \mathbb{P}(T_z \in [n - n^{9/10}, n])$$
$$\leq \mathbb{P}(\tilde{T}_z - \tilde{T}_z^l > n^{9/10}) + \mathbb{P}(T_z \in [n - n^{9/10}, n])$$
$$\leq Kn^{-1/30} \log n + Kn^{-1/30} \log^{10/3} n \leq Kn^{-1/30} \log^{10/3} n$$

by Lemmas 4.1(b) and 4.2(a). Part (b) is done in the same way, but we use Lemmas 4.1(a) and 4.2(c). □

The next two preparatory lemmas give estimates similar to those in Lemmas 4.1 and 4.2, and will be used to prove Lemma 4.6, where we show that for points $y$ and $z$ that are neither too close to the boundary of $C_n(z)$ nor to each other, it is unlikely that $y \notin \tilde{C}_n(z)$ if $y \in C_n(z)$.

We introduce stopping times reminiscent of those defined at the beginning of this section but which are concerned with the time at which points $y$ and $z$ first lie (or "almost" lie) in different components of the complement of the Brownian or random walk path. If $y, z \in \mathbb{C}$,

$$T_{y,z} = \inf\{t \geq 0 : \exists 0 \leq s \leq t \text{ with } B(s) = B(t),$$
$$|\arg_y(B(s), B(t))| \neq |\arg_z(B(s), B(t))|\},$$
$$\tilde{T}_{y,z} = \inf\{t \geq 0 : \exists 0 \leq r \leq t \text{ with } S(r) = S(t),$$
$$|\arg_y(S(r), S(t))| \neq |\arg_z(S(r), S(t))|\},$$
$$T_{y,z}^l = T_{y,z}^l(n) = \inf\{t \geq 0 : \exists 0 \leq s \leq t \text{ with } d(B(s), B(t)) \leq 3n^{1/4} \log^2(n),$$
$$|\arg_y(B(s), B(t)) - \arg_z(B(s), B(t))| \geq 3\pi/2\},$$
$$\tilde{T}_{y,z}^l = \tilde{T}_{y,z}^l(n) = \inf\{t \geq 0 : \exists 0 \leq r \leq t \text{ with } d(S(r), S(t)) \leq 3n^{1/4} \log^2(n),$$
$$|\arg_y(B(s), B(t)) - \arg_z(B(s), B(t))| \geq 3\pi/2\}.$$

We now give the analogues of Lemmas 4.1 and 4.2 for the quantities we just defined. As the arguments are the same, we only indicate which modifications from the proofs of Lemmas 4.1 and 4.2 are needed.

LEMMA 4.4. *There exists a constant $K > 0$ such that for any $y, z \in \mathbb{Z}^2$, any $n \geq 1$ and any $a > 1/2$:*

(a) $\mathbb{P}(T_{y,z} - T_{y,z}^l > n^a) \leq K(n^{(1-2a)/24} \log n + |y - z|^{-1/6} n^{1/24} \log n).$
(b) $\mathbb{P}(\tilde{T}_{y,z} - \tilde{T}_{y,z}^l > n^a) \leq K(n^{(1-2a)/24} \log n + |y - z|^{-1/6} n^{1/24} \log n).$

PROOF. The proof is essentially the same as that of Lemma 4.1. The only difference is that we need to take into account the relative positions of $y$ and $x^l_{y,z}$, as well as those of $z$ and $x^l_{y,z}$. We define $m(y,z,a) = \min\{|z - x^l_{y,z}|, |y - x^l_{y,z}|, n^{a/2}/\log n\}$. By symmetry, we can assume $|z - x^l_{y,z}| \leq |y - x^l_{y,z}|$. We then find that if $m(y,z,a) \geq 100n^{1/4}\log^2 n$, then $\mathbb{P}(T_{y,z} - T^l_{y,z} > n^a) \leq Kn^{1/8}(m(y,z,a))^{-1/2}\log n$ and if $m(y,z,a) = |z - x^l_{y,z}|$, then

$$\mathbb{P}(T_{y,z} - T^l_{y,z} > n^a) \leq K\left(\frac{|z - x^l_{y,z}|}{\min\{|y - x^l_{y,z}|, n^{a/2}/\log n\}}\right)^{1/4}.$$

Looking separately at the cases $|y - x^l_{y,z}| \leq n^{a/2}/\log n$ and $n^{a/2}/\log n \leq |y - x^l_{y,z}|$, we find that

$$\mathbb{P}(T_{y,z} - T^l_{y,z} > n^a) \leq K(n^{(1-2a)/24}\log n + |y - x^l_{y,z}|^{-1/6}n^{1/24}\log n).$$

Since we assumed that $|z - x^l_{y,z}| \leq |y - x^l_{y,z}|$, we have $|y - x^l_{y,z}| \geq |y - z|/2$, which concludes the proof. □

LEMMA 4.5. *There exists a constant $K > 0$ such that for any $y, z \in \mathbb{Z}^2, n \in \mathbb{N}$ and $a > 1/2$:*

(a) $\mathbb{P}(T_{y,z} \in [n - n^a, n]) \leq Kn^{(a-1)/3}\log^{11/3} n$,
(b) $\mathbb{P}(T^l_{y,z} \in [n - n^a, n]) \leq Kn^{(a-1)/3}\log^{11/3} n$,
(c) $\mathbb{P}(\tilde{T}_{y,z} \in [n - n^a, n]) \leq Kn^{(a-1)/3}\log^{11/3} n$,
(d) $\mathbb{P}(\tilde{T}^l_{y,z} \in [n - n^a, n]) \leq Kn^{(a-1)/3}\log^{11/3} n$.

PROOF. The proof is virtually the same as that of Lemma 4.2 except the definitions of $\Gamma$ and $\Gamma^l$ must now be $\Gamma = \Gamma_{y,z,w} = \inf\{t \geq \Phi_w : |B(t) - w| = \min\{|z - w|, |y - w|\}\}$, and $\Gamma^l = \Gamma^l_{y,z,w} = \sup\{t \leq T^l_z : |B(t) - w| = \min\{|z - w|, |y - w|\}\}$. □

Recall the definitions made in (4.8).

LEMMA 4.6. *For any given $y, z \in \mathbb{Z}^2$:*

(a) $\mathbb{P}(\{y \in C_n(z) \setminus \tilde{C}_n(z)\}; \mathcal{B}(z); \mathcal{B}(y); \mathcal{P})$
$\quad \leq K\max\{n^{-1/30}\log^2 n, |y - z|^{-1/6}n^{1/24}\log n\}$.
(b) $\mathbb{P}(\{y \in \tilde{C}_n(z) \setminus C_n(z)\}; \tilde{\mathcal{B}}(z); \tilde{\mathcal{B}}(y); \mathcal{P})$
$\quad \leq K\max\{n^{-1/30}\log^2 n, |y - z|^{-1/6}n^{1/24}\log n\}$.

PROOF. (a) Note that $y$ and $z$ lie in different components of $\mathbb{C} \setminus B[0,t]$ if and only if $t \geq T_{y,z}$. Exactly as in Lemma 4.3, $\{\mathcal{B}(z); \mathcal{B}(y); \mathcal{P}\} \subset \{T^l_{y,z} \leq$

$\tilde{T}_{y,z}\}$, and so, for $a > 1/2$,

$$\mathbb{P}(\{y \in C_n(z) \setminus \tilde{C}_n(z)\}; \mathcal{B}(z); \mathcal{B}(y); \mathcal{P})$$
$$\leq \mathbb{P}(y \in C_n(z); T^l_{y,z} < n - n^a) + \mathbb{P}(\tilde{T}_{y,z} \in [n - n^a, n]).$$

The result now follows from Lemmas 4.4(a) and 4.5(c).

The proof of (b) is the same but uses Lemmas 4.4(b) and 4.5(a). □

We are now ready to attack the core of the argument. The following proposition shows that under the coupling, if $C_n(z)$ and $\tilde{C}_n(z)$ are both finite, one of them is "large enough," and $z$ is not too close to their boundaries, then the difference between the areas of $C_n(z)$ and $\tilde{C}_n(z)$ is usually small, relatively to these areas. We define

$$\Delta(z) = ||C_n(z)| - |\tilde{C}_n(z)||,$$

the difference in area between the Brownian motion and random walk hole containing a given point $z \in \mathbb{Z}^2$ and

$$\mathcal{E}(z) = \{|C_n(z)| < \infty; |\tilde{C}_n(z)| < \infty\},$$

the condition that both the Brownian and random walk components containing $z$ are finite. Recall that $\mathcal{P} = \{\sup_{0 \leq t \leq n} |B(t) - S(2t)| \leq n^{1/4} \log^2 n\}$ and $\mathcal{B}(z) = \{d(z, \partial C_n(z)) \geq 100 n^{1/4} \log^2 n\}$.

PROPOSITION 4.1. *For every $K_1 > 0$, there exists a constant $K_2 > 0$ such that for any $z \in \mathbb{C}$, $n \in \mathbb{N}$ and $0 < b < 1/30$,*

$$\mathbb{P}(\{\Delta(z) \geq K_1 n^{1-b}\}; \mathcal{B}(z); \mathcal{E}(z); \mathcal{P}) \leq K_2 n^{b-1/30} \log^4 n.$$

PROOF. We will prove the proposition under the additional assumption that

$$\mathcal{N} = \left\{\sup_{0 \leq t \leq n} |B(t)| \leq \sqrt{n} \log n\right\}$$

is satisfied, which will suffice since $\mathbb{P}(\mathcal{N}^c)$ decays rapidly by Lemma 2.1. We will use the abbreviation

$$\mathcal{G} = \mathcal{G}(z, n) := \mathcal{B}(z) \cap \mathcal{E}(z) \cap \mathcal{P} \cap \mathcal{N}$$

and show that

$$\mathbb{E}[\Delta(z); \mathcal{G}] \leq K n^{1-1/30} \log^4 n.$$

Once we have this, the proposition follows from Chebyshev's inequality. Note that

(4.11)
$$\begin{aligned}\Delta(z) \leq &\#\{y \in \mathbb{Z}^2 : \operatorname{Sq}(y) \cap \partial C_n(z) \neq \varnothing\} \\ &+ \#\{y \in \mathbb{Z}^2 : \operatorname{Sq}(y) \cap \partial \tilde{C}_n(z) \neq \varnothing\} \\ &+ \#\{y \in \mathbb{Z}^2 : y \in C_n(z) \setminus \tilde{C}_n(z)\} \\ &+ \#\{y \in \mathbb{Z}^2 : y \in \tilde{C}_n(z) \setminus C_n(z)\}.\end{aligned}$$

The first two terms on the right-hand side are related to the Hausdorff dimension of the Brownian frontier and were dealt with in Section 3. If we modify (4.11) slightly and "thicken" the boundary, it follows from Proposition 3.1 that

(4.12)
$$\begin{aligned}\mathbb{E}[\Delta(z); \mathcal{G}] \leq &\mathbb{E}[\#\{y \in \mathbb{Z}^2 : d(y, \partial C_n(z)) \leq 100 n^{1/4} \log^2 n\}] \\ &+ \mathbb{E}[\#\{y \in \mathbb{Z}^2 : d(y, \partial \tilde{C}_n(z)) \leq 100 n^{1/4} \log^2 n\}] \\ &+ \mathbb{E}\left[\sum_{y \in \mathbb{Z}^2} \mathbb{1}\{y \in C_n(z) \setminus \tilde{C}_n(z); \mathcal{G}; \mathcal{B}(y)\} \right. \\ &\left. \qquad + \mathbb{1}\{y \in \tilde{C}_n(z) \setminus C_n(z); \mathcal{G}; \tilde{\mathcal{B}}(y)\}\right] \\ \leq &K n^{11/12}(\log n)^{17/3} + \sum_{y \in \mathbb{Z}^2} \mathbb{P}(\{y \in C_n(z) \setminus \tilde{C}_n(z)\}; \mathcal{G}; \mathcal{B}(y)) \\ &+ \sum_{y \in \mathbb{Z}^2} \mathbb{P}(\{y \in \tilde{C}_n(z) \setminus C_n(z)\}; \mathcal{G}; \tilde{\mathcal{B}}(y)).\end{aligned}$$

By Lemma 4.6(a) and the definition of $\mathcal{G}$,

$$\begin{aligned}\sum_{y \in \mathbb{Z}^2} &\mathbb{P}\{y \in C_n(z) \setminus \tilde{C}_n(z); \mathcal{G}; \mathcal{B}(y)\} \\ &\leq \sum_{y \in D(\sqrt{n} \log n) \cap \mathbb{Z}^2} \mathbb{P}(\{y \in C_n(z) \setminus \tilde{C}_n(z)\}; \mathcal{B}(z); \mathcal{B}(y); \mathcal{P}) \\ &\leq K \sum_{|y| \leq n^{1/4}} 1 + \sum_{n^{1/4} \leq |y| \leq n^{9/20}} n^{1/24} |y|^{-1/6} \log n \\ &\qquad + \sum_{n^{9/20} \leq |y| \leq n^{1/2} \log n} n^{-1/30} \log^2 n \\ &\leq K n^{1 - 1/30} \log^4 n,\end{aligned}$$

and the same bound holds for the second sum in (4.12). Thus,

$$\mathbb{E}[\Delta(z); \mathcal{G}] \leq K n^{11/12} (\log n)^{17/3} + K n^{1-1/30} \log^4 n \leq K n^{1-1/30} \log^4 n. \quad \square$$

## 5. Main results.

PROPOSITION 5.1. *There exists a probability space $(\Omega, \mathcal{F}, \mathbb{P})$ containing a planar simple random walk $S$ and a planar standard Brownian motion $B$ such that, if for $\delta \in \mathbb{R}$,*

$N_n(\delta) = \#\{\text{connected components of } \mathbb{C} \setminus B[0,n] \text{ of area } \geq n^{1-\delta}\},$

$\tilde{N}_n(\delta) = \#\{\text{connected components of } \mathbb{C} \setminus S[0,2n] \text{ of area } \geq n^{1-\delta}\},$

$\hat{N}_n(\delta) = \#\{\text{connected components of } \mathbb{Z}^2 \setminus S[0,2n] \text{ of cardinality } \geq n^{1-\delta}\},$

*then for every $\varepsilon > 0$, every $0 < \delta < 1/60$,*

(5.1) $\qquad \mathbb{P}(|\tilde{N}_n(\delta) - N_n(\delta)| > \varepsilon N_n(\delta)) \to 0 \qquad \text{as } n \to \infty$

*and*

(5.2) $\qquad \mathbb{P}(|\hat{N}_n(\delta) - N_n(\delta)| > \varepsilon N_n(\delta)) \to 0 \qquad \text{as } n \to \infty.$

REMARK. We will only prove (5.1), as the proof of (5.2) is virtually the same. The reason why the same proof works is that what would cause the random walk hole containing a point $z$ to have an area that is substantially different from the lattice area of the random walk lattice hole containing that point is their boundary behavior and that boundary effect is eliminated in Proposition 4.1 using the work done in Section 3.

Before proving (5.1), we introduce the notation used in this section, as well as the extensions of [10] which are essential to our proof. Although most of the quantities we are about to define depend on $\delta$, we will usually let this dependence be implicit to keep the notation as light as possible:

$N_{[a,b)} = N_{[a,b)}(n) = \#\{\text{conn. components of } \mathbb{C} \setminus B[0,n] \text{ with area in } [a,b)\},$

$\tilde{N}_{[a,b)} = \#\{\text{conn. components of } \mathbb{C} \setminus S[0,2n] \text{ with area in } [a,b)\},$

$I_n = I_n(\delta) = [n^{1-\delta}, \infty), \qquad N_n(\delta) = N_{I_n}, \qquad \tilde{N}_n(\delta) = \tilde{N}_{I_n}.$

We let

$$c = c(\varepsilon) = 1 + \varepsilon$$

and define for all $j \geq -1$,

$$I_{j,n}^R = [n^{1-\delta} c^{j+1} (1+\varepsilon^2)^{-1}, n^{1-\delta} c^{j+1}),$$

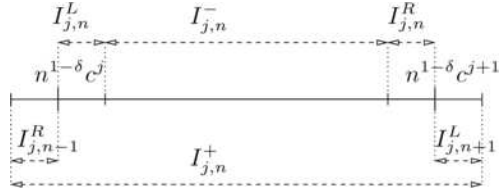

Fig. 5. *The intervals needed in the proof of Proposition* 5.1.

for all $j \geq 0$,

$$I_{j,n} = [n^{1-\delta}c^j, n^{1-\delta}c^{j+1}), \qquad I_{j,n}^L = [n^{1-\delta}c^j, n^{1-\delta}c^j(1+\varepsilon^2)),$$
$$I_{j,n}^- = I_{j,n} \setminus (I_{j,n}^L \cup I_{j,n}^R), \qquad I_{j,n}^+ = I_{j,n} \cup I_{j-1}^R \cup I_{j+1}^L.$$

See Figures 5 and 6. The number of components in the corresponding intervals will be

(5.3)
$$N_{j,n} = N_{I_{j,n}}, \qquad \tilde{N}_{j,n} = \tilde{N}_{I_{j,n}},$$
$$N_{j,n}^L = N_{I_{j,n}^L}, \qquad N_{j,n}^R = N_{I_{j,n}^R}, \qquad N_{j,n}^\pm = N_{I_{j,n}^\pm}.$$

We let

$Z_{j,n}^\pm$ = the set of all components of $\mathbb{C} \setminus B[0,n]$ with area in $I_{j,n}^\pm$,

$\tilde{Z}_{j,n}$ = the set of all components of $\mathbb{C} \setminus S[0,2n]$ with area in $I_{j,n}$,

$Z_{j,n}$ = the set of all components of $\mathbb{C} \setminus B[0,n]$ with area in $I_{j,n}$.

The details on the results which we state now and which follow from [10], can be found in the Appendix. If we let

$$\gamma_n = \gamma_n(\delta) = \frac{2\pi n^\delta}{\log^2(n^\delta)},$$

then for every $K > 0$,

(5.4) $\qquad \mathbb{P}(|N_n(\delta) - \gamma_n| \geq K\gamma_n) \to 0 \qquad$ as $n \to \infty$.

In particular, for every $\delta > 0, K_1 < 1 < K_2$,

(5.5) $\quad \lim_{n \to \infty} \mathbb{P}(N_n(\delta) < K_1 \gamma_n) = 0 \quad$ and $\quad \lim_{n \to \infty} \mathbb{P}(N_n(\delta) > K_2 \gamma_n) = 0$.

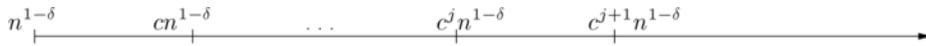

Fig. 6. *Splitting $[n^{1-\delta}, \infty)$ into small finite intervals.*

We define

$$\gamma_{j,n} = \frac{2\pi n^\delta \log c}{c^j \log^2 c^j/(\pi n^\delta)}, \qquad \gamma_{j,n}^{LR} = \frac{2\pi n^\delta \log(1+\varepsilon^2)}{c^j \log^2 c^j/(\pi n^\delta)},$$

(5.6)

$$\gamma_{j,n}^{-} = \frac{2\pi n^\delta \log c/(1+\varepsilon^2)^2}{c^j \log^2 c^j/(\pi n^\delta)},$$

and point out that for every $\delta > 0$, there exists a $K = K(\delta) > 0$ such that for all $n \geq 1, \varepsilon > 0$ and $0 \leq j \leq [\frac{\delta \log n}{2 \log c}]$,

(5.7) $$|\gamma_{j,n}^{LR} - \varepsilon \gamma_{j,n}^{-}| \leq K\varepsilon^2 \gamma_{j,n}^{-}.$$

Then (A.5) and (A.6) imply that for every $K_1 > 0$ and $\delta > 0$, there is a constant $K_2 = K_2(K_1, \delta) > 0$ such that for all $n$ large enough, $\varepsilon > 0$ small enough, and all $0 \leq j \leq [\frac{\delta \log n}{2 \log c}]$,

(5.8) $$\mathbb{P}(|N_{j,n} - \gamma_{j,n}| \geq K_1 \gamma_{j,n}) \leq K_2 \log^{-3/2} n,$$

and the same inequality holds if we replace the pair $(N_{j,n}, \gamma_{j,n})$ by $(N_{j,n}^L, \gamma_{j,n}^{LR})$, $(N_{j,n}^R, \gamma_{j,n}^{LR})$, or $(N_{j,n}^{-}, \gamma_{j,n}^{-})$. We are now ready to prove the proposition:

PROOF OF PROPOSITION 5.1. We will show that for all $0 < \delta < 1/60$ and $\varepsilon > 0$,

$$\mathbb{P}(|\tilde{N}_n(\delta) - N_n(\delta)| > 11\varepsilon N_n(\delta)) \to 0 \qquad \text{as } n \to \infty.$$

We first note that "most" of the Brownian holes of area greater than $n^{1-\delta}$ have an area which is "close to" $n^{1-\delta}$. For instance, as can be seen from (5.4), the number of holes for $B$ with area in the interval $[n^{1-\delta}, n^{1-\delta/2}]$ is typically of greater order of magnitude than the number of Brownian holes with area in the interval $[n^{1-\delta/2}, \infty)$. Recall that $c = 1 + \varepsilon$. If we let $m = m(n, \varepsilon, \delta) = [\frac{\delta \log n}{2 \log c}]$, then $n^{1-\delta} c^m \leq n^{1-\delta/2} \leq n^{1-\delta} c^{m+1}$ and so

(5.9)
$$\mathbb{P}(|\tilde{N}_n(\delta) - N_n(\delta)| > 11\varepsilon N_n(\delta))$$
$$\leq \mathbb{P}(|\tilde{N}_{[n^{1-\delta}, n^{1-\delta}c^m)} - N_{[n^{1-\delta}, n^{1-\delta}c^m)}| > 10\varepsilon N_n(\delta))$$
$$+ \mathbb{P}(|\tilde{N}_{[n^{1-\delta}c^m, \infty)} - N_{[n^{1-\delta}c^m, \infty)}| > \varepsilon N_n(\delta)).$$

The last term of (5.9) goes to 0 as $n$ goes to $\infty$. Indeed, if $|\tilde{N}_{[n^{1-\delta}c^m, \infty)} - N_{[n^{1-\delta}c^m, \infty)}| > \varepsilon N_n(\delta)$, then either $\tilde{N}_{[n^{1-\delta}c^m, \infty)} > \varepsilon N_n(\delta)$ or $N_{[n^{1-\delta}c^m, \infty)} > \varepsilon N_n(\delta)$. We also know from (5.5) that $\mathbb{P}(N_n(\delta) \leq \frac{2n^\delta}{\log^2 n^\delta}) = o(1)$. By observing that if $N_{[n^{1-\delta}c^m, \infty)} > \varepsilon \frac{2n^\delta}{\log^2 n^\delta}$, the total area enclosed by $B$ is greater

than $\varepsilon \frac{2n^{1+\delta/2}}{\log^2 n^\delta}$, and that the same holds for $S$, we can conclude that

$$\mathbb{P}(|\tilde{N}_{[n^{1-\delta}c^m,\infty)} - N_{[n^{1-\delta}c^m,\infty)}| > \varepsilon N_n(\delta))$$
$$\leq \mathbb{P}\left(\sup_{0\leq t\leq n} |B(t)| \geq \varepsilon \frac{2n^{(1+\delta/2)/2}}{\sqrt{\pi}\log n^\delta}\right)$$
$$+ \mathbb{P}\left(\sup_{0\leq t\leq n} |S(2t)| \geq \varepsilon \frac{2n^{(1+\delta/2)/2}}{\sqrt{\pi}\log n^\delta}\right) + o(1) \to 0 \qquad \text{as } n \to \infty,$$

by (2.4) and (2.5). It now suffices to show that

$$\lim_{n\to\infty} \mathbb{P}(|\tilde{N}_{[n^{1-\delta},n^{1-\delta}c^m)} - N_{[n^{1-\delta},n^{1-\delta}c^m)}| > 10\varepsilon N_n(\delta)) = 0.$$

Observe that

$$\mathbb{P}(|\tilde{N}_{[n^{1-\delta},n^{1-\delta}c^m)} - N_{[n^{1-\delta},n^{1-\delta}c^m)}| > 10\varepsilon N_n(\delta))$$
$$\leq \mathbb{P}\left(\sum_{j=0}^{m-1} |\tilde{N}_{j,n} - N_{j,n}| > 10\varepsilon \sum_{j=0}^{m-1} N_{j,n}\right)$$
$$\leq \sum_{j=0}^{m-1} \mathbb{P}(|\tilde{N}_{j,n} - N_{j,n}| > 10\varepsilon N_{j,n}).$$

Because of the definition of $m$, it now suffices to prove that for every $\varepsilon > 0, 0 < \delta < \frac{1}{60}$, there is a $K > 0$ such that for all $0 \leq j \leq m-1$,

(5.10) $$\mathbb{P}(|\tilde{N}_{j,n} - N_{j,n}| > 10\varepsilon N_{j,n}) \leq K\psi(n),$$

where $\psi(n) = o(\log^{-1} n)$. Recall that we defined

$\tilde{Z}_{j,n} = \{$connected components of $\mathbb{C} \setminus S[0,2n]$ with area in $I_{j,n}\}$,

$Z_{j,n} = \{$connected components of $\mathbb{C} \setminus B[0,n]$ with area in $I_{j,n}\}$.

Unfortunately, we cannot use Proposition 4.1 to show (5.10) quite yet. The problem is that the fact that $|\tilde{N}_{j,n} - N_{j,n}| > 10\varepsilon N_{j,n}$ does not imply anything about $||Z_{j,n}| - |\tilde{Z}_{j,n}||$ and so estimates about the latter cannot be used to show (5.10). To make things work, we need the interval containing the areas of random walk holes to be strictly included in the interval containing the areas of Brownian motion holes (or vice versa). The key in creating such a situation is to observe that

(5.11)
$$\mathbb{P}(|\tilde{N}_{j,n} - N_{j,n}| > 10\varepsilon N_{j,n})$$
$$\leq \mathbb{P}(\tilde{N}_{j,n} \leq (1-2\varepsilon)N_{j,n}^-) + \mathbb{P}(\tilde{N}_{j,n} \geq (1+2\varepsilon)N_{j,n}^+)$$
$$+ \mathbb{P}(N_{j,n}^L \geq 4\varepsilon N_{j,n}^-) + \mathbb{P}(N_{j,n}^R \geq 4\varepsilon N_{j,n}^-)$$
$$+ \mathbb{P}(N_{j+1}^L \geq 7/2\varepsilon N_{j,n}) + \mathbb{P}(N_{j-1}^R \geq 7/2\varepsilon N_{j,n}).$$

In words, if $\tilde{N}_{j,n}$ is much greater than $N_{j,n}$, then either it is somewhat greater than $N_{j,n}^+$, or there are many Brownian holes with area in $I_{j,n}^+ \setminus I_{j,n}$. This scheme will work because we have defined $I_{j,n}^+$ in such a way that the Lebesgue measure of $I_{j,n}^+ \setminus I_{j,n}$ is of a smaller order of magnitude than the Lebesgue measure of $I_{j,n}$.

$$\mathbb{P}\{N_{j,n}^L \geq 4\varepsilon N_{j,n}^-\}$$
$$\leq \mathbb{P}(N_{j,n}^L \geq 4\varepsilon N_{j,n}^-; N_{j,n}^L \leq \tfrac{3}{2}\gamma_{j,n}^{LR}; N_{j,n}^- \geq \tfrac{1}{2}\gamma_{j,n}^-)$$
$$+ \mathbb{P}(N_{j,n}^L \geq \tfrac{3}{2}\gamma_{j,n}^{LR}) + \mathbb{P}(N_{j,n}^- \leq \tfrac{1}{2}\gamma_{j,n}^-)$$
$$\leq \mathbb{1}\{\tfrac{3}{2}\gamma_{j,n}^{LR} \geq 2\varepsilon\gamma_{j,n}^-\} + \mathbb{P}(N_{j,n}^L \geq \tfrac{3}{2}\gamma_{j,n}^{LR}) + \mathbb{P}(N_{j,n}^- \leq \tfrac{1}{2}\gamma_{j,n}^-).$$

We know from (5.7) that if $\varepsilon$ is small enough, $\tfrac{3}{2}\gamma_{j,n}^{LR} \leq 2\varepsilon\gamma_{j,n}^-$, so the first term is 0. The second and the third are $\mathcal{O}((\log n)^{-3/2})$, uniformly for $0 \leq j \leq m-1$, by (5.8), and the last three terms of (5.11) can be bounded in the same way. Therefore,

$$\mathbb{P}(|\tilde{N}_{j,n} - N_{j,n}| > 10\varepsilon N_{j,n})$$
(5.12)
$$\leq \mathbb{P}(\tilde{N}_{j,n} \leq (1-2\varepsilon)N_{j,n}^-)$$
$$+ \mathbb{P}(\tilde{N}_{j,n} \geq (1+2\varepsilon)N_{j,n}^+) + \mathcal{O}((\log n)^{-3/2}).$$

We define

$$\Delta_j^+ = |\tilde{Z}_{j,n}| - |Z_{j,n}^+|, \qquad \Delta_j^- = |Z_{j,n}^-| - |\tilde{Z}_{j,n}|.$$

To bound the first term on the right-hand side of (5.12), note that

$$|Z_{j,n}^-| \geq N_{j,n}^- c^j n^{1-\delta}(1+\varepsilon^2) \quad \text{and} \quad |\tilde{Z}_{j,n}| \leq \tilde{N}_{j,n} c^{j+1} n^{1-\delta},$$

so that if $\tilde{N}_{j,n} \leq (1-2\varepsilon)N_{j,n}^-$, then

$$\Delta_j^- \geq N_{j,n}^- n^{1-\delta} c^j [1 + \varepsilon^2 - (1-2\varepsilon)c] = N_{j,n}^- n^{1-\delta} c^j (\varepsilon + 3\varepsilon^2).$$

Since $N_{j,n}^- \in \mathbb{N}$ and we know from (5.8) that $\mathbb{P}(N_{j,n}^- = 0) \leq K(\log n)^{-3/2}$, we see that

(5.13) $\quad \mathbb{P}(\tilde{N}_{j,n} \leq (1-2\varepsilon)N_{j,n}^-; \Delta_j^- \leq \varepsilon n^{1-\delta}) \leq K(\log n)^{-3/2}.$

The obvious inequality $|A| - |B| \leq |A \setminus B|$ for sets $A, B \in \mathbb{C}$ and (5.13) imply that

$$\mathbb{P}(\tilde{N}_{j,n} \leq (1-2\varepsilon)N_{j,n}^-) \leq \mathbb{P}(|Z_{j,n}^- \setminus \tilde{Z}_{j,n}| \geq \varepsilon n^{1-\delta}) + \mathcal{O}((\log n)^{-3/2}).$$

But

$$Z_{j,n}^- \setminus \tilde{Z}_{j,n} \subset \bigcup \mathrm{Sq}(z),$$

where the union is over the union of $\{z \in \mathbb{Z}^2 : d(z, \partial Z_{j,n}^-) \leq 100 n^{1/4} \log^2 n\}$ and $\{z \in \mathbb{Z}^2 : |C_n(z)| \in I_{j,n}^-; |\tilde{C}_n(z)| \notin I_{j,n}; d(z, \partial Z_{j,n}^-) \geq 100 n^{1/4} \log^2 n\}$. We can now combine the different pieces of our work and write

$$\mathbb{E}[|Z_{j,n}^- \setminus \tilde{Z}_{j,n}|; \mathcal{P}; \mathcal{N}]$$

$$\leq \mathbb{E}\left[\sum_{|z| \leq \sqrt{n} \log n} \mathbb{1}\{d(z, \partial Z_{j,n}^-) \leq 100 n^{1/4} \log^2 n\}\right]$$

$$+ \sum_{|z| \leq \sqrt{n} \log n} \mathbb{P}(\{|C_n(z)| \in I_{j,n}^-; |\tilde{C}_n(z)| \notin I_{j,n}\}; \mathcal{B}(z); \mathcal{P})$$

$$\leq K n^{(33+\delta)/36} \log^{14/3} n + \sum_{|z| \leq \sqrt{n} \log n} \mathbb{P}(\{\Delta(z) \geq \varepsilon^2 n^{1-\delta}\}; \mathcal{B}(z); \mathcal{E}(z); \mathcal{P})$$

$$+ \sum_{|z| \leq \sqrt{n} \log n} \mathbb{P}(\{|C_n(z)| \in I_{j,n}^-; |\tilde{C}_n(z)| = \infty\}; \mathcal{B}(z); \mathcal{P})$$

$$\leq K(n^{11/12+\delta/36} \log^{14/3} n + n^{1+\delta-1/30} \log^6 n + n^{1-1/30} \log^{16/3} n)$$

$$\leq K n^{1+\delta-1/30} \log^6 n,$$

by Propositions 3.2 and 4.1 and Lemma 4.3. The sums are over elements of $\mathbb{Z}^2$ and the constant $K$ may depend on $\varepsilon$. Therefore,

$$\mathbb{P}(|Z_{j,n}^- \setminus \tilde{Z}_{j,n}| \geq \varepsilon n^{1-\delta})$$

$$\leq \mathbb{P}(\{|Z_{j,n}^- \setminus \tilde{Z}_{j,n}| \geq \varepsilon n^{1-\delta}\}; \mathcal{P}; \mathcal{N}) + \mathbb{P}(\mathcal{P}^c) + \mathbb{P}(\mathcal{N}^c)$$

$$\leq K n^{1+\delta-1/30} \log^6 n / \varepsilon n^{1-\delta} + K n^{1-b \log n} + K n^{-\log n/2}$$

$$\leq (K/\varepsilon) n^{2\delta-1/30} \log^6 n.$$

For every $\varepsilon > 0$ and $0 < \delta < 1/60$, this decays faster than $(\log n)^{-3/2}$. The second term of (5.12) is bounded in the same way. It now suffices to look back at (5.10) to see that the proof is complete. $\square$

Given this proposition, the proof of Theorem 1.1 is straightforward:

$$\mathbb{P}(|\tilde{N}_n(\delta) - 2\pi \gamma_n| > \varepsilon \gamma_n)$$

$$\leq \mathbb{P}\left(|\tilde{N}_n(\delta) - N_n(\delta)| > \frac{\varepsilon}{2} \gamma_n\right) + \mathbb{P}\left(|N_n(\delta) - 2\pi \gamma_n| > \frac{\varepsilon}{2} \gamma_n\right)$$

$$\leq \mathbb{P}\left(|\tilde{N}_n(\delta) - N_n(\delta)| > \frac{\varepsilon}{4} N_n(\delta)\right) + \mathbb{P}(\gamma_n < N_n(\delta)/2)$$

$$+ \mathbb{P}\left(|N_n(\delta) - 2\pi \gamma_n| > \frac{\varepsilon}{2} \gamma_n\right).$$

From Proposition 5.1, (5.5) and (5.4), we know that for every $\varepsilon > 0$, every $0 < \delta < 1/60$, each of the 3 terms goes to 0 as $n$ goes to infinity. The same holds if we replace $\tilde{N}_n(\delta)$ by $\hat{N}$, which concludes the proof of the theorem.

REMARK. The statement of Theorem 1.1 can be written in the following equivalent form: If $A(n) = n^{1-\delta}, 0 < \delta \leq \delta_0$, then with the notation defined in (1.2),

$$(5.14) \qquad \frac{\log^2(n/A(n))}{n/A(n)} \mathcal{H}_n(A(n)) \xrightarrow{P} 2\pi \qquad \text{as } n \to \infty,$$

$$(5.15) \qquad \frac{\log^2(n/A(n))}{n/A(n)} \mathcal{L}_n(A(n)) \xrightarrow{P} 2\pi \qquad \text{as } n \to \infty.$$

A closer look at the proof shows that in fact, (5.14) and (5.15) hold for arbitrary $A(n)$, as long as $n^{1-\delta_0} \leq A(n)$ and $A(n) = o(n)$. A straightforward consequence of this and Lemma 2.1 is that if $n^{1-\delta_0} \leq A(n)$ and $A(n) = o(n)$, then

$$\log \mathbb{E}[\mathcal{H}_n(A(n))] \sim \log(n/A(n))$$

and

$$\log \mathbb{E}[\mathcal{L}_n(A(n))] \sim \log(n/A(n)).$$

**6. Small holes.** Our interest in this problem was spurred by the observation by Mandelbrot (see [11]), based on computer simulations, that while the number of random walk holes and the number of large lattice holes are governed by a power function with exponent 2, the exponent for small lattice holes is 5/3.

The remark at the end of Section 5 shows that if $A(n) \geq n^{1-\delta}$ and $A(n) = o(n)$, then $\log \mathbb{E}[\mathcal{L}_n(A(n))] \sim \log(nA(n)^{-\xi/2})$, with $\xi = 2$. The number $\xi$ is what Mandelbrot calls the exponent.

A formal definition can be made as follows: Given a sequence of intervals $\{I_n\}_{n \in \mathbb{N}}$, a sequence of functions $\{f_n\}: I_n \to \mathbb{R}$ has *exponent* $\xi$ over the intervals $J_n \subset I_n$ if for any two sequences $\{A(n)\}_{n \in \mathbb{N}}, \{A'(n)\}_{n \in \mathbb{N}}$ with $A(n), A'(n) \in J_n$,

$$(6.1) \qquad \log(f_n(A(n))/f_n(A'(n))) \sim -\frac{\xi}{2} \log(A(n)/A'(n)).$$

With the notation above, Mandelbrot's observation can be described as follows:

There exist sequences $A_0, A_1$, and $A_2$ with $\log A_2(n) \sim \log n$, such that:

- $\mathbb{E}[\mathcal{L}_n(\cdot)]$ has exponent 5/3 over $[1, A_0(n)]$.
- $\mathbb{E}[\mathcal{L}_n(\cdot)]$ has exponent strictly between 5/3 and 2 over $[A_0(n), A_1(n)]$.

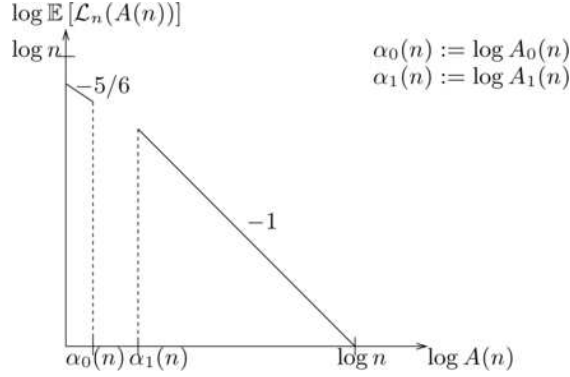

FIG. 7. *Two regimes. The log–log plot yields a slope of $-1$ for large holes and $-5/6$ for small holes, corresponding to exponents of $2$ and $5/3$, respectively.*

- $\mathbb{E}[\mathcal{L}_n(\cdot)]$ has exponent 2 over $[A_1(n), A_2(n)]$.

As mentioned above, Theorem 1.1 implies that there is a $\delta_0 > 0$ such that the third point is true with $A_1(n) = n^{1-\delta_0}$ and any $A_2(n)$ satisfying $A_2(n) = o(n)$, $\log A_2(n) \sim \log n$, and the same holds if $\mathbb{E}[\mathcal{L}_n(\cdot)]$ is replaced by $\mathbb{E}[\mathcal{H}_n(\cdot)]$.

It turns out we can say much more about $A_1$ if the situation is exactly as described above in the three points. Assuming that situation, it is easy to see that if $A_1(n) = n^\varepsilon$ for some $\varepsilon > 0$, then $\lim_{n\to\infty} \log \mathbb{E}[\mathcal{L}_n(1)]/\log n < 1$ [see Figure 7 and recall that $\mathcal{L}_n(1)$ is the number of lattice holes of lattice area $\geq 1$]. We will see in the lemma below that in fact, $\log \mathbb{E}[\mathcal{L}_n(1)] \sim \log n$, thus showing that in the situation described above, $A_1(n)$ would necessarily grow more slowly than any power function. This would imply that $\log A_1(n) = o(\log n)$, so that the region in Figure 7 over which the exponent $5/3$ may be observed would become a smaller and smaller fraction of the domain of $\log \mathbb{E}[\mathcal{L}_n(\cdot)]$, as $n \to \infty$. The same is true if we replace $\mathbb{E}[\mathcal{L}_n(\cdot)]$ by $\mathbb{E}[\mathcal{H}_n(\cdot)]$.

We devote the remainder of this section to proving that $\log \mathbb{E}[\mathcal{L}_n(1)] \sim \log n$ and $\log \mathbb{E}[\mathcal{H}_n(1)] \sim \log n$, and conclude by stating a conjecture based on Mandelbrot's picture and the results derived in this paper.

We say that two (lattice) holes made by $S[0, 2n]$, $H$ and $H'$, are equivalent if there exists a $z \in \mathbb{Z}^2$ such that $H' = H + z$. We choose a representative $H$ for each equivalence class by requesting that the first point of $\partial H$ in the lexicographic order be the origin and call the set of these representatives $\mathcal{R}_1$ for holes and $\mathcal{R}_2$ for lattice holes. Given a representative $H$ and a point $z \in \mathbb{Z}^2$, we define $H_z := H + z$. We will write $\bar{H}_z = H_z \cup \partial H_z$ and $\tau_z(H) = \inf\{k \geq 0 : S(k) \in \bar{H}_z\}$. For each $H$ in $\mathcal{R}_1$ or $\mathcal{R}_2$, we define $N_H = \#\{z \in \mathbb{Z}^2 : H_z \text{ is a (lattice) hole made by } S[0, 2n]\}$.

LEMMA 6.1. *For any given hole $H \in \mathcal{R}_1$ or lattice hole $H \in \mathcal{R}_2$, there exists a $K > 0$ such that for all $n \geq 1$,*

$$\mathbb{E}[N_H] \geq K \frac{n}{\log^2 n}.$$

PROOF. We give the proof for $H \in \mathcal{R}_1$. The case where $H \in \mathcal{R}_2$ is done in the same way. Suppose the boundary of $H$ is composed of $2k$ line segments of length 1 and that $n \geq 8k$. Then $\mathbb{E}[N_H]$ is equal to

$$\sum_{z \in \mathbb{Z}^2} \mathbb{P}(H_z \text{ is a hole for } S[0, 2n])$$

$$= \sum_{z \in \mathbb{Z}^2} \sum_{l=0}^{2n-2k} \mathbb{P}(H_z \text{ is a hole for } S[0, 2n]; \tau_z(H) = l)$$

$$\geq \sum_{z \in \mathbb{Z}^2} \sum_{l=[n/2]}^{[3n/2]} \mathbb{P}(S[0, l-1] \cap \bar{H}_z = \varnothing;$$

$$S[l, l+2k] = \partial H_z; S[l+2k, 2n] \cap H_z = \varnothing)$$

$$\geq K \min_{w \in \partial H \cap \mathbb{Z}^2} \mathbb{P}^w(S[1, [3n/2]] \cap \bar{H} = \varnothing) n$$

$$\times \min_{w \in \partial H \cap \mathbb{Z}^2} \mathbb{P}^w(S[0, [3n/2]] \cap H = \varnothing),$$

where the last inequality is obtained from the Markov property, by considering the time reversal of $S[0, l-1]$, and by translation invariance of simple random walk. The fact that for any hole $H$ and $w \notin H$, there is a $K > 0$ such that $\mathbb{P}^w(S[0, n] \cap H = \varnothing) \sim K/\log n$, which follows from Lemma 2.3.1 in [3], concludes the proof. □

COROLLARY 6.1.

$$\log \mathbb{E}[\mathcal{H}_n(1)] \sim \log n \quad \text{and} \quad \log \mathbb{E}[\mathcal{L}_n(1)] \sim \log n.$$

PROOF. This follows immediately from Lemma 6.1 and the fact that there exists a constant $K > 0$ such that $\mathbb{E}[\mathcal{H}_n(1)] \leq Kn \log^4 n$ and $\mathbb{E}[\mathcal{L}_n(1)] \leq Kn \log^4 n$, which is a direct consequence of Lemma 2.1. □

We now summarize what we believe should be the global picture, based on the work presented in this paper and our understanding of Mandelbrot's observation: The exponent for the expected number of (lattice) holes is 2, except for the case of small lattice holes, for which the exponent is 5/3. The exponent 5/3, however, only holds for lattice holes whose lattice areas grow more slowly than any power function.

CONJECTURE. *There exist nondecreasing sequences $A_0, A_1$ and $A_2$, with $1 < A_0(n) < A_1(n) < A_2(n)$, where $A_1(n) = o(\log n)$ and $\log A_2(n) \sim \log n$, such that*

$$\log\left(\frac{\mathbb{E}[\mathcal{H}_n(A(n))]}{\mathbb{E}[\mathcal{H}_n(A'(n))]}\right) \sim -\frac{\xi}{2}\log\left(\frac{A(n)}{A'(n)}\right)$$

*and*

$$\log\left(\frac{\mathbb{E}[\mathcal{L}_n(A(n))]}{\mathbb{E}[\mathcal{L}_n(A'(n))]}\right) \sim -\frac{\eta}{2}\left(\frac{A(n)}{A'(n)}\right),$$

*with:*

- $\xi = 2$ if $1 \leq A(n), A'(n) \leq A_2(n)$,
- $\eta = 2$ if $A_1(n) \leq A(n), A'(n) \leq A_2(n)$,
- $\eta = \frac{5}{3}$ if $1 \leq A(n), A'(n) \leq A_0(n)$.

## APPENDIX: FROM SMALL TO LARGE BROWNIAN HOLES

In [10], estimates for the expectation and the variance of $|U_\eta|$, the area of $U_\eta$, are derived, where

$$U_\eta = \{y \in \mathbb{C} : \pi(\lambda\eta)^2 \leq |C_1(y)| \leq \pi\eta^2\},$$

$\lambda < 1$, and $C_1(y)$ is defined as in (2.2). In particular, it is shown that the variance is of smaller order of magnitude than the second moment. The two estimates which are relevant to us are the following:

1.

(A.2) $$\mathbb{E}[|U_\eta|] = \frac{\pi|\log \lambda|}{|\log \eta|^2}\left(1 + \mathcal{O}\left(\frac{\log|\log \eta|}{|\log \eta|^{1/2}}\right)\right),$$

where $\mathcal{O}(\cdot)$ is for $\eta \to 0$, but the implied constant may depend on $\lambda$.

2. There exists a constant $K > 0$ such that for every $\eta \in (0, 1/4)$,

(A.3) $$\mathrm{Var}(|U_\eta|) \leq K|\log \eta|^{-11/2}.$$

If we let $A_\zeta = \{y \in \mathbb{C} : \zeta n \leq |C_n(y)| \leq c\zeta n\}$, where $c = 1 + \varepsilon > 1$, then the scaling property of Brownian motion allow us to deduce the following from (A.2) and (A.3):

$$\mathbb{E}[|A_\zeta|] = \frac{2\pi n \log c}{\log^2(c\zeta/\pi)}\left(1 + \mathcal{O}\left(\frac{\log|\log(c\zeta)|}{|\log(c\zeta)|^{1/2}}\right)\right),$$

where $\mathcal{O}(\cdot)$ is for $\zeta \to 0$, but the implied constant may depend on $c$,

(A.4) $$\mathrm{Var}(|A_\zeta|) \leq Kn^2|\log(c\zeta)|^{-11/2},$$

for all $\zeta$ with $c\zeta \in (0, \pi/16)$.

We can now easily translate these facts into results about the number of components of area lying in a certain interval, rather than the total area covered by these components. This just requires dividing the total area by the area of a single component. Since for $y \in A_\zeta$, $|C_n(y)|$ can take any value in $[\zeta n, c\zeta n]$, we have an additional error term.

If $\delta > 0$ and $m = [\frac{\delta \log n}{2 \log c}]$, then for every $j \leq m - 1$, components of area in $I_{j,n}$ [see (1.5)] have area less than $n^{1-\delta/2} \leq \frac{\pi}{16} n$, so we can use (A.4). With the definitions of (5.3) and (5.6), we have for $0 \leq j \leq m - 1$,

$$\text{(A.5)} \qquad \mathbb{E}[N_{j,n}] = \gamma_{j,n}\left(1 - r\varepsilon + \mathbb{O}\left(\frac{\log|\log(c^j n^{-\delta})|}{|\log(c^j n^{-\delta})|^{1/2}}\right)\right),$$

$$\text{(A.6)} \qquad \text{Var}(N_{j,n}) \leq K\gamma_{j,n}^2 |\log(c^j n^{-\delta})|^{-3/2},$$

$$\mathbb{E}[N_{j,n}^-] = \gamma_{j,n}^-\left(1 + r\varepsilon + \mathbb{O}\left(\frac{\log|\log(c^j n^{-\delta})|}{|\log(c^j n^{-\delta})|^{1/2}}\right)\right),$$

$$\text{Var}(N_{j,n}^-) \leq K(\gamma_{j,n}^-)^2 |\log(c^j n^{-\delta})|^{-3/2},$$

$$\mathbb{E}[N_{j,n}^L] = \gamma_{j,n}^{LR}\left(1 + r\varepsilon^2 + \mathbb{O}\left(\frac{\log|\log(c^j n^{-\delta})|}{|\log(c^j n^{-\delta})|^{1/2}}\right)\right),$$

$$\text{Var}(N_{j,n}^L) \leq K(\gamma_{j,n}^{LR})^2 |\log(c^j n^{-\delta})|^{-3/2},$$

and for $-1 \leq j \leq m - 1$,

$$\mathbb{E}[N_{j,n}^R] = \gamma_{j,n}^{LR}\left(1 + r\varepsilon^2 + \mathbb{O}\left(\frac{\log|\log(c^j n^{-\delta})|}{|\log(c^j n^{-\delta})|^{1/2}}\right)\right),$$

$$\text{Var}(N_{j,n}^R) \leq K(\gamma_{j,n}^{LR})^2 |\log(c^j n^{-\delta})|^{-3/2},$$

where $|r| < 2$ and $K$ and the constants of $\mathcal{O}$ may depend on $\varepsilon$. $\mathcal{O}(\cdot)$ is for $n \to \infty$. (5.8) and its analogues follow directly from the set of equations above and (5.4) is a consequence of a scaled version of the main result in [10]: If $N(u)$ is the number of connected components of $\mathbb{C} \setminus B[0,1]$ of area greater than $u$, then

$$\lim_{u \to 0} u(\log u)^2 N(u) = 2\pi \qquad \text{a.s.}$$

**Acknowledgments.** The author wishes to express his sincere gratitude to Greg Lawler for his guidance and support, to Chris Burdzy for directing the author to [12] and [10], as well as to the reviewer for many helpful comments.

TUFTS UNIVERSITY
503 BOSTON AVE
MEDFORD, MASSACHUSETTS 02155
USA
E-MAIL: christian.benes@tufts.edu